\documentclass[a4paper]{amsart}

\usepackage{amscd,amsthm,amsfonts,amssymb,amsmath}
\usepackage[colorlinks=true]{hyperref}

    \newcommand{\BA}{{\mathbb {A}}} 
    \newcommand{\BC}{{\mathbb {C}}} 
     \newcommand{\BF}{{\mathbb {F}}}
    \newcommand{\BG}{{\mathbb {G}}} \newcommand{\BH}{{\mathbb {H}}}

     \newcommand{\BP}{{\mathbb {P}}}
    \newcommand{\BQ}{{\mathbb {Q}}} \newcommand{\BR}{{\mathbb {R}}}

     \newcommand{\BZ}{{\mathbb {Z}}}

    \newcommand{\CA}{{\mathcal {A}}}

    \newcommand{\CG}{{\mathcal {G}}} \newcommand{\CH}{{\mathcal {H}}}

    \newcommand{\CO}{{\mathcal {O}}}

     \newcommand{\RN}{{\mathrm {N}}}

    \newcommand{\ad}{{\mathrm{ad}}}
    \newcommand{\alg}{{\mathrm{alg}}}
    \newcommand{\Aut}{{\mathrm{Aut}}}

    \newcommand{\End}{{\mathrm{End}}}

    \newcommand{\Gal}{{\mathrm{Gal}}} \newcommand{\GL}{{\mathrm{GL}}}
    
    \newcommand{\Hom}{{\mathrm{Hom}}}
    
    \renewcommand{\Im}{{\mathrm{Im}}}

    \newcommand{\Ker}{{\mathrm{Ker}}}

    \newcommand{\ord}{{\mathrm{ord}}} \newcommand{\rank}{{\mathrm{rank}}}

    \renewcommand{\mod}{\ \mathrm{mod}\ }

    \newcommand{\SL}{{\mathrm{SL}}}

    \newcommand{\tor}{{\mathrm{tor}}}
    \newcommand{\RTr}{{\mathrm{Tr}}}
    
    \newcommand{\Vol}{{\mathrm{Vol}}}

\newcommand{\matrixx}[4]{\begin{pmatrix}
#1 & #2 \\ #3 & #4
\end{pmatrix} }        

    \font\cyr=wncyr10

    \newcommand{\Sha}{\hbox{\cyr X}}\newcommand{\wt}{\widetilde}
    \newcommand{\wh}{\widehat}
    
    \newcommand{\pair}[1]{\langle {#1} \rangle}

    \newcommand{\ov}{\overline}

    \newcommand{\sk}{\medskip}
    \newcommand{\lra}{\longrightarrow}
    \newcommand{\ra}{\rightarrow} 
    \newcommand{\lto}{\longmapsto}\newcommand{\bs}{\backslash}
    \newcommand{\nequiv}{\equiv\hspace{-9.5pt}/\ }
    \newcommand{\s}{\sk\noindent}

    \theoremstyle{plain}
    \newtheorem{thm}{Theorem}[section] 
    \newtheorem{lem}[thm]{Lemma}  \newtheorem{prop}[thm]{Proposition}
     \newtheorem{defn}[thm]{Definition}

\theoremstyle{remark} \newtheorem{remark}[thm]{Remark}
\theoremstyle{remark} 
\theoremstyle{remark} 

    \newcommand{\adeles}{ad\'{e}les~}

    \numberwithin{equation}{section}

\begin{document}

\title{Congruent Numbers and Heegner Points}

\author{Ye Tian}
\address{Academy of Mathematics and Systems Science,
Morningside center of Mathematics, Chinese Academy of Sciences,
Beijing 100190} \email{ytian@math.ac.cn}
\thanks{The author was supported by NSFC grant 11031004, 973 Program
2013CB834202, and The Chinese Academy of Sciences The Hundred
Talents Program.}

\maketitle

\tableofcontents

\section{Introduction and Main Results}

A positive integer is called a {\em congruent number} if it is the
area of a right-angled  triangle, all of whose sides have rational
length. The problem of determining which positive integers are
congruent is buried in antiquity (see Chapter 9 of Dickson
\cite{Dickson}), with it long being known that the numbers 5, 6, and
7 are congruent. Fermat proved that 1 is not a congruent number, and
similar arguments show that also 2 and 3 are not congruent numbers.
No algorithm has ever been proven for infallibly deciding whether a
given integer $n\geq 1$ is congruent. The reason for this is that it
can easily been that an integer $n\geq 1$ is congruent if and only
if there exists a point $(x, y)$, with $x$ and $y$ rational numbers
and $y\neq 0$, on the elliptic curve $ny^2 = x^3-x$.  Moreover,
assuming $n$ to be square free, a classical calculation of root
numbers shows that the complex L-function of this curve has  zero of
odd order at the center of its critical strip precisely when $n$
lies in one of the residue classes of $5$, $6$, and $7$ modulo  $8$.
Thus, in particular, the unproven conjecture of Birch and
Swinnerton-Dyer predicts that every positive integer lying in the
residue classes of $5$, $6$, and $7$ modulo $8$ should be a
congruent number. The aim of this paper is to prove the following
partial results in this direction.

\begin{thm}\label{Easy Theorem}For any given integer $k\geq 0$, there are infinitely
many square-free  congruent  numbers  with exactly $k+1$ odd prime
divisors in each residue class of $5, 6$, and $7$ modulo $8$.
\end{thm}

\begin{remark}\label{remark1}
 The above result when $k=0$ is due to Heegner \cite{Heegner}, Birch
\cite{Birch},  Stephens \cite{Stephens}, and completed by Monsky
\cite{Monsky},  and that when $k=1$ is due to Monsky \cite{Monsky}
and Gross \cite{Tunnell}. Actually Heegner is the first
mathematician who found (in \cite{Heegner}) a method to construct
fairly general solutions to cubic Diophantine equations. The method
of this paper is based on his construction.
\end{remark}

In addition to Theorem \ref{Easy Theorem}, we have the following
result on the conjecture of Birch and Swinnerton-Dyer. For any
abelian group $A$ and an integer $d\geq 1$, we write $A[d]$ for the
kernel of multiplication by $d$ on $A$.

\begin{thm}\label{Main Theorem 2}Let $k\geq 0$ be an integer and
$n=p_0p_1\cdots p_k$ a product of distinct odd primes  with
$p_i\equiv 1\mod 8$ for $1\leq i\leq k$.  Assume  that the ideal
class group $\CA$ of the field $K=\BQ(\sqrt{-2n})$ satisfies the
condition:
\begin{equation}
\dim_{\BF_2}(\CA[4]/\CA[2])=
\begin{cases}0, \quad
&\text{if $n\equiv \pm 3\mod 8$},\\
1, &\text{otherwise.} \end{cases}
\end{equation}
Let $m=n$ or $2n$ such that $m \equiv 5, 6$,  or $7\mod 8$. Let
$E^{(m)}$ be the elliptic curve $my^2=x^3-x$ over $\BQ$. Then we
have
$$\rank_\BZ E^{(m)}(\BQ)=1=\ord_{s=1}L(E^{(m)}, s).$$
Moreover,  the Shafarevich-Tate group of $E^{(m)}$ is finite and has
odd cardinality.
\end{thm}

\begin{remark}
The work of Perrin-Riou \cite{Perrin-Riou} and Kobayshi
\cite{Kobayshi} shows that the order of the $p$-primary subgroup of
the Tate-Shafarevich group of $E^{(m)}$ is as predicted by the
conjecture of Birch and Swinnerton-Dyer for all primes $p$ with $(p,
2m)=1$. At present, it is unknown whether the same statement holds
for the primes $p$ dividing $2m$, so that the full
Birch-Swinnerton-Dyer conjecture is still not quite completely known
for the curves $E^{(m)}$. However, toward to the conjecture for
$p=2$ we have  Theorem \ref{Main Theorem 3} below in viewing of
Gross-Zagier formula.
\end{remark}

The condition (1.1) on $\CA[4]/\CA[2]$ in Theorem \ref{Main Theorem
2} allows us to complete the first 2-descent and to show that the
2-Selmer group of $E^{(m)}$ modulo the 2-torsion subgroup of
$E^{(m)}$  is $\BZ/2\BZ$ (see Lemma 5.1). It follows that
$$\rank_\BZ E^{(m)}(\BQ) +\dim_{\BF_2} \Sha(E^{(m)}/\BQ)[2]=1,$$
and therefore that $\rank_\BZ E^{(m)}(\BQ)$ is either 0 or 1. Any
one of the parity conjecture for Mordell-Weil group and the
finiteness conjecture for Shafarevich-Tate group predicts that the
elliptic curve $E^{(m)}$ has Mordell-Weil group of rank 1. Therefore
the BSD conjecture predicts that the analytic rank is 1. Then the
generalization of Gross-Zagier formula  predicts that the height of
a Heegner divisor is non-zero, so this Heegner divisor class should
have infinite order. However, we shall follow a different path, and,
always assuming (1.1), we shall prove independently of any
conjectures that this Heegner divisor does indeed have infinite
order, and the Tate-Shafarevich group is finite of odd order. Our
method uses induction on the number of primes dividing the congruent
number $m$, Kolyvagin's Euler system, and a generalization of the
Gross-Zagier formula.

Let $E$ be the elliptic curve $y^2=x^3-x$ so that $E^{(m)}$ is a
quadratic twist of $E$.  It is well known that the only rational
torsion on $E^{(m)}$ is the subgroup $E^{(m)}[2]$ of 2-torsion. Let
$E(\BQ(\sqrt{m}))^-$ denote the subgroup of those points in
$E(\BQ(\sqrt{m}))$ which are mapped to their negative by the
non-trivial element of the Galois group of $\BQ(\sqrt{m})$ over
$\BQ$. Then the map which sends $(x, y)$ to $(x, \sqrt{m}y)$ defines
an isomorphism form $E^{(m)}(\BQ)$ onto $E(\BQ(\sqrt{m}))^-$. Thus
$m$ will be congruent if and only if we can show that
$E(\BQ(\sqrt{m}))^-$ is strictly larger than $E[2]$. Note that
$E^{(m)}$ and $E^{(-m)}$ are isomorphic over $\BQ$.

The modular curve $X_0(32)$  of level $\Gamma_0(32)$ has genus 1 and
is defined over $\BQ$.  Its associated Riemann surface structure  is
given by the complex uniformization
$$X_0(32)(\BC)=\Gamma_0(32)\bs (\CH \cup \BP^1(\BQ)),$$  where $\CH$
is the upper half complex plane, and we write $[z]$ for the point on
the curve defined by any $z\in \CH\cup \BP^1(\BQ)$. It is easy to
see that $[\infty]$ is defined over $\BQ$.  The elliptic curve $E$
has conductor 32 and there is a degree 2 modular parametrization $f:
X_0(32)\ra E$ mapping $[\infty]$ to $0$. Such $f$ is unique up to
multiplication by $-1$ because the elliptic curve $(X_0(32),
[\infty])$ has only one rational torsion point of order $2$ (see
Proposition \ref{Main Prop}). Let $n=p_0p_1\cdots p_k$ and $m$ be
integers as in Theorem \ref{Main Theorem 2}. Let $K=\BQ(\sqrt{-2n})$
and $H$ its Hilbert class field. Let $m^*=(-1)^{\frac{n-1}{2}}m$ and
 $\chi$ the abelian character over $K$ defining the unramified extension
$K(\sqrt{m^*})$. Define the point $P\in X_0(32)$ to be
$[i\sqrt{2n}/8]$ if $n\equiv 5\mod 8$, and to be
$[(i\sqrt{2n}+2)/8]$ if $n\equiv 6$ or $7\mod 8$. Both Theorem
\ref{Easy Theorem} and Theorem \ref{Main Theorem 2} will follow from
the following main theorem of the paper.
\begin{thm}\label{Main Theorem 3}
Let $n=p_0p_1\cdots p_k$ and $m$ be integers as in Theorem \ref{Main
Theorem 2}. Then the point $f(P)\in E$ is defined over $H(i)$; and
the $\chi$-component of $f(P)$, defined by
$$P^\chi(f):=\sum_{\sigma \in \Gal(H(i)/K)} f(P)^\sigma
\chi(\sigma),$$satisfies
$$P^\chi(f)\in 2^{k+1} E(\BQ(\sqrt{m^*}))^- \quad\text{and}\quad
P^\chi(f)\notin 2^{k+2}E(\BQ(\sqrt{m^*}))^-+E[2].$$In particular,
$P^\chi(f)\in E(\BQ(\sqrt{m^*}))^-\cong E^{(m)}(\BQ)$ is of infinite
order and $m$ is a congruent number.
\end{thm}

We now explain our method in the case $n\equiv 5\mod 8$ in details.
Other cases are similar. Let $p_0\equiv 5\mod 8$ and $p_i\equiv
1\mod 8$, $1\leq i\leq k$,  be distinct primes. Let $n=p_0p_1\cdots
p_k$ and $K=\BQ(\sqrt{-2n})$. The theory of complex multiplication
implies that
$$z:=f(P)+(1+\sqrt{2}, 2+\sqrt{2})$$ is a point on $E$ defined over
the Hilbert class field $H$ of $K$, even though neither $f(P)$ nor
the 4-torsion point $(1+\sqrt{2}, 2+\sqrt{2})$ on $E$ is defined
over $H$. Note that if we use $-f$ instead of $f$, then we still
obtain a $H$-rational point $-f(P)+(1+\sqrt{2}, 2+\sqrt{2})=-z+(1,
0)$. The desired Heegner point is defined by taking trace of $z$
from $H$ to $K(\sqrt{n})$
\begin{eqnarray}y_n:=\RTr_{H/K(\sqrt{n})}z\in
E(K(\sqrt{n})).\end{eqnarray} It turns out that $y_n$ is actually
defined over $\BQ(\sqrt{n})$. Moreover, $y_n$ (resp. $2y_n$) belongs
to $E(\BQ(\sqrt{n}))^-$ if $k\geq 1$ (resp. $k=0$). Now it is easy
to see that the point $P^\chi(f)$, defined in Theorem \ref{Main
Theorem 3}, is equal to $4y_n$.

The condition (1.1) in Theorem \ref{Main Theorem 2}, in the case
$p_0\equiv 5\mod 8$,  is equivalent to that the Galois group
$\Gal(H/H_0)\cong 2\CA$ has odd cardinality where $H_0=K(\sqrt{p_0},
\cdots, \sqrt{p_k})$ is the genus field of $K$. We will show that
the point $y_n$ is of infinite order for $n=p_0p_1\cdots p_k$
satisfying the condition (1.1) in Theorem \ref{Main Theorem 2}. When
$k=0$, classical arguments show that $y_n$ is of infinite order
(see, for example, \cite{Monsky}). However, when $k$ is at least 1,
we will prove by induction on $k$ that, provided $n=p_0p_1\cdots
p_k$ satisfies condition (1.1) in Theorem \ref{Main Theorem 2}, the
point $y_n$ belongs to $2^{k-1} E(\BQ(\sqrt{n}))^-+E[2]$, but does
not belong to $2^k E(\BQ(\sqrt{n}))^-+E[2]$. This clearly shows that
$y_n$ must be of infinite order. (Note that condition (1.1) holds
automatically when $k=0$). We now give some more details on how
these arguments are carried through in detail.

In fact, we find a relation of $y_n$ with other Heegner divisors.
Now assume that $k\geq 1$ and let $y_0=\RTr_{H/H_0}z\in E(H_0)$. It
turns out that $y_0\in E(H_0^+)$ where $H_0^+=H_0\cap \BR$. For any
positive divisor $d$ of $n$ divisible by $p_0$, let
$y_d=\RTr_{H/K(\sqrt{d})} z$, which actually belongs to
$E(\BQ(\sqrt{d}))^-$. The points $y_d$'s with $p_0|d|n$ and $y_0$
are related by the following relation:
\begin{equation}\sum_{p_0|d|n} y_d=\begin{cases}
2^k y_0, \qquad &\text{if $k\geq 2$},\\
2^k y_0+\#2\CA \cdot (0, 0), &\text{if
$k=1$}.\end{cases}\end{equation}

In the next, for any proper divisor $d$ of $n$ divisible by $p_0$,
we need to know the 2-divisibility of $y_d$ in the Mordell-Weil
group $E(\BQ(\sqrt{d}))^-$. To do this, we  similarly construct a
point $y_d^0\in E(\BQ(\sqrt{d}))$ with $K$ replaced by
$K_0=\BQ(\sqrt{-2d})$, whose 2-divisibility is understood by
induction hypothesis. We can reduce the comparison of
2-divisibilities of $y_d$ and $y_d^0$ to the comparison of their
heights via Kolyvagin's result. The heights of these two points are
related to central derivative L-values via Gross-Zagier formula
Theorem 1.2 in \cite{YZZ}.  The comparison of heights of these two
points is further reduced to the comparison of two central L-values,
which is given by Zhao in \cite{Zhao1}. It turns out from the
comparison and induction hypothesis that $y_d\in 2^k
E(\BQ(\sqrt{d}))^-+E[2]$ for all proper divisors $d$ of $n$. It
follows from the equality (1.3) that
$$y_n=2^k \left(y_0-\sum_{p_0|d|n, d\neq n} y_d'\right) +t$$
for some $y_d'\in E(\BQ(\sqrt{d}))^-$ and $t\in E[2]$. It can then
be shown by additional arguments (see the proof of Theorem \ref{Main
Theorem A}) that this implies that  $y_n \in
2^{k-1}E(\BQ(\sqrt{n}))^-+E[2]$. The fact that $y_n\notin 2^k
E(\BQ(\sqrt{n}))^-+E[2]$ with $n$ satisfying the condition (1.1)
follows from the same algebraic ingredient as in the initial case
$k=0$ and some ramification argument. Note that $4y_n=P^\chi(f)$ and
then Theorem \ref{Main Theorem 3} follows in the case $n\equiv 5\mod
8$.

\begin{remark}\label{remark3}By a conjecture of Goldfeld \cite{Goldfeld} or Katz-Sarnak
\cite{Katz-Sarnak},  combined with Coates-Wiles' result
\cite{Coates-Wiles}, almost all positive integers $n\equiv 1,
2,3\mod 8$ are non-congruent numbers.

It is known (\cite{Feng} and \cite{Li-Tian}) that for any given
integer $k\geq 0$, there are infinitely many square-free
non-congruent positive integers with exactly $k+1$ odd prime
divisors in each residue class of $1$, $2$, and $3$ modulo $8$.  In
fact,  let $n=p_0p_1\cdots p_k$ be a product of distinct odd primes
with $p_i\equiv 1\mod 8$ for $1\leq i\leq k$ satisfying the
condition (1.1) in Theorem \ref{Main Theorem 2}.  Let $m=n$ or $2n$
such that $m\equiv 1, 2, 3\mod 8$. Then $m$ is non-congruent if
$p_0\nequiv 1\mod 8$.  Moreover, if $p_0\equiv 1\mod 8$, then $n$ is
non-congruent provided the additional assumption $\prod_{i=0}^k
\left(\left(\frac{2}{p_i}\right)_4\cdot (-1)^{(p_i-1)/8}\right)=-1$.

The above non-congruent numbers are constructed easily by minimizing
the 2-Selmer groups attached to 2-isogenies of $E^{(m)}$  and taking
2-part of the Shafarevich-Tate group into account in the case
$p_0\equiv 1\mod 8$.
\end{remark}

\s{\bf Notations and Conventions}. We often work with the imaginary
quadratic field $K=\BQ(\sqrt{-2n})$ where $n$ is a square-free
positive odd integer. We fix an embedding of the algebraic closure
$\bar{K}$ of $K$ in $\BC$. Let $\CO_K$ denote the ring of integers
in $K$. Then $\CO_K=\BZ+\BZ w$ with $w=i\sqrt{2n}$. Let $K^{ab}$
denote the maximal abelian extension of $K$. Let $\wh{\BZ}=\prod_p
\BZ_p$  and $\wh{K}=K\otimes_\BZ \wh{\BZ}$ the finite \adeles of
$K$. Denote by
$$[\quad , K^{ab}/K]:\quad \wh{K}^\times/K^\times \lra \Gal(K^{ab}/K),
$$the Artin reciprocity law, and similarly for $[\quad, \BQ^{ab}/\BQ]$.
We also often write $\sigma_t=[t, K^{ab}/K]$. For each prime $p|2n$,
let $\varpi_p$ be a uniformizer $\sqrt{-2n}$ of the local field
$K_p$ of $K$ at the unique prime above $p$;  for each $0<d|n$, let
$\varpi_d=\prod_{p|d} \varpi_p\in \wh{K}^\times$.  We often use the
convention $\varpi=\varpi_2\in K_2^\times$. For $t\in \wh{K}^\times$
(resp. $\wh{\BQ}^\times$), we denote by $t_2\in K_2^\times$  (resp.
$\BQ_2$) its component at $2$.

We will also need Gauss' genus theory for the imaginary quadratic
field $K$. Denote by $H$ the Hilbert class field of $K$ and $\CA$
the ideal class group of $K$. Suppose that $n$ has exact $k+1$ prime
divisors: $n=\prod_{j=0}^kp_j$.  Let $H_0=K(\sqrt{p_0^*},
\sqrt{p_1^*}, \cdots, \sqrt{p_k^*})\subset H$ be its genus field of
$K$, where $p^*=(-1)^{(p-1)/2} p$ so that $p^*\equiv 1\mod 4$.
Sometime we identify the ideal class group $\CA$ of $K$ with the
class group $\wh{K}^\times/K^\times \wh{\CO}_K^\times$. Consider the
exact sequence
$$0\lra \CA[2]\lra \CA \stackrel{\times 2}{\lra} \CA \lra \CA/2\CA \lra 0.$$
Gauss' genus theory says the following
\begin{itemize}
\item[(i)] the subgroup $\CA[2]$ consists of ideal classes of $(d, \sqrt{-2n})$,
corresponding to classes of $\varpi_d$ in
$\wh{K}^\times/K^\times\wh{\CO}_K^\times$, where $d$ runs over all
positive divisors of $n$. Therefore $\CA[2]$ has cardinality
$2^{k+1}$.
\item[(ii)]  under the class field theory isomorphism
$\sigma: \CA \simeq \Gal(H/K)$, the subgroup $2\CA \simeq
\Gal(H/H_0)$, i.e.  the class of $t\in 2\CA$  if and only if
$\sigma_t$ fixes all $\sqrt{p_j^*}, 0\leq j\leq k$.
\end{itemize}
By Gauss'  quadratic reciprocity law, for each $0<d|n$ and $p|n$,
$\sigma_{\varpi_d}$ fixes $\sqrt{p^*}$ iff
$\left(\frac{d}{p}\right)=1$ for $p\nmid d$ and
$\left(\frac{2n/d}{p}\right)=1$ for $p|d$.

It is easy to see that $i\notin H$ and therefore the restriction map
gives  natural isomorphisms $\Gal(H(i)/K(i))\cong \Gal(H/K)$. Let
$\CO_2$ denote the ring of integers in the local field $K_2$. Its
unit group $\CO_2^\times$ is generated by $-1, 5, 1+\varpi$ as a
$\BZ_2$-module. We often need the structure of the Galois group
$\Gal(H(i)/\BQ)\cong \Gal(H(i)/K)\rtimes \{1, c\}$, where $c$ is the
complex conjugation and
$$\Gal(H(i)/K)\cong \wh{K}^\times \big/K^\times\wh{\CO}_K^{\times (2)}
U_2, \quad \text{with}\quad U_2=\BZ_2^\times (1+2\varpi \CO_2).$$
Here the supscript in $\wh{\CO}_K^{\times (2)}$ means the component
above 2 is removed. Note that $U_2\subset \CO_2^\times$ is generated
by $-1, 5, (1+\varpi)^2$. Thus the Galois group $\Gal(H(i)/H)$ is
generated by $\sigma_{1+\varpi}$. The group $\Gal(H(i)/\BQ)$ is
generated by $\Gal(H(i)/H_0(i)) \cong 2\CA$, the complex conjugation
$c$, and elements representing $\Gal(H_0(i)/K)$. For example, when
$2\CA\cap \CA[2]=0$, $\Gal(H_0(i)/K)$ is represented by
$\sigma_{1+\varpi}$ and $\sigma_{\varpi_d}$ with all $0<d|n$.

\s{\bf Acknowledgment}. The author thanks John Coates, Xinyi Yuan,
Shouwu Zhang, Wei Zhang, and the referee for many useful discussions
and comments. In the original version of this paper, the author used
the construction of Heegner points given by employing the modular
parametrization of $E$ via the modular curve $X(8)$, following the
original work of Heegner and Monsky. The current simpler
construction using the curve $X_0(32)$ arose out of discussions with
Xinyi Yuan.

The author thanks Keqin Feng, Delang Li, Mingwei Xu, and Chunlai
Zhao for bringing him to this beautiful topic when he was in a
master degree program. The author thanks John Coates,  Benedict
Gross, Victor Kolyvagin, Yuan Wang, Lo Yang, Shing-Tung Yau, and
Shouwu Zhang  for constant encouragement during preparation of this
work.

\section{Modular Parametrization and CM Points}

Let $E$ be the elliptic curve with Weierstrass equation $y^2=x^3-x$.
It is known that $E$ has conductor $32$.  Let $f: X_0(32)\lra E$ be
a fixed modular parametrization over $\BQ$ of degree $2$ mapping the
cusp $[\infty]$ at the infinity on $X_0(32)$ to the zero  element
$0\in E$. In this section, we will construct suitable CM points on
$E$ associated to the modular parametrization $f$ and the imaginary
quadratic field $K=\BQ(\sqrt{-2n})$ with $n$ a positive square-free
odd integer. We will show that $E':=(X_0(32), [\infty])$ is an
elliptic curve with Weierstrass equation $y^2=x^3+4x$. Before give
construction of points on $E$, we need set up the correspondence of
torsion points of $E'$ over $\BQ(i)$ between their $(x,
y)$-coordinates and their modular expressions.

We now recall the following standard notation. Let $\CH$ be the
upper half complex plane, on which the subgroup $\GL_2^+(\BR)$ of
elements of $\GL_2(\BR)$ with positive determinant acts by linear
fractional transformations. Let $\Gamma_0(32)$ be the subgroup of
$\SL_2(\BZ)$ consisting of all matrices $\matrixx{a}{b}{c}{d}$ with
$c\equiv 0\mod 32$, which  acts on $\CH \cup \BP^1(\BQ$ by linear
fractional transformation.  Denote by $Y_0(32)$ the modular curve of
level $\Gamma_0(32)$ over $\BQ$ and $X_0(32)$ its projective closure
over $\BQ$. Then the underlying compact Riemann surface of $X_0(32)$
is given as:
$$X_0(32)(\BC)=Y_0(32)(\BC) \cup S,$$where
$$Y_0(32)(\BC)=\Gamma_0(32)\bs \CH,\qquad S=\Gamma_0(32)\bs \BP^1(\BQ).$$
For each $z\in \CH \cup \BP^1(\BQ)$, let $[z]$ denote the point on
$X_0(32)(\BC)$ represented by $z$. The set $S$ consists of 8 cusps:
$$[\infty],\quad [0], \quad [-1/2],\quad [-1/16],\quad  [-1/4],\quad [-3/4],\quad
[1/8], \quad [-1/8],$$ where the first 4 cusps are defined over
$\BQ$ and the later 4 ones have the field of definition $\BQ(i)$.
The curve $X_0(32)$ has genus one and thus we have an elliptic curve
$E':=(X_0(32), [\infty])$ over $\BQ$ with the cusp $[\infty]$  as
its zero element.

\begin{prop}\label{Add 1} The elliptic curve $E'=(X_0(32), \infty)$
has complex multiplication by $\BZ[i]$ and  Weierstrass equation
$y^2=x^3+4x$ such that the cusp $[0]=(2, 4)$ in $(x,
y)$-coordinates.  Moreover, the set $S$ of cusps on $X_0(32)$ is
exactly $E'[(1+i)^3]$.
\end{prop}
\begin{proof}
Define $N$ to be the normalizer of $\Gamma_0(32)$ in $\GL_2^+(\BR)$
and let $Z(\BR)$ denote the center of $\GL_2^+(\BR)$. Let
$\Aut(X_0(32)(\BC))$ denote the group of automorphisms of
$X_0(32)(\BC)$ and $\Aut(X_0(32)(\BC), S)$ its subgroup of
automorphisms $t$ satisfying $t(S)=S$. Then the action of $N\subset
\GL_2^+(\BR)$ on $\CH\cup \BP^1(\BQ)$ induces a homomorphism
$$T: N\lra \Aut(X_0(32)(\BC), S)$$with kernel
$Z(\BR)\Gamma_0(32)$.

Now, as is very well known, every element of $\Aut(X_0(32)(\BC))$ is
of form
$$t_{\alpha, \epsilon}(x)=\epsilon(x)+\alpha,$$
where $\epsilon$ belongs to the group $\Aut(E'_\BC)$ of
automorphisms of the elliptic curve $E'_\BC$ (i.e. ones with
$\epsilon([\infty])=[\infty]$), and $\alpha$ is some point in
$E'(\BC)$.

Now consider the matrices
$$A=\matrixx{0}{1}{-32}{0}, \quad B=\matrixx{1}{1/4}{0}{1}, \quad
C=AB^2=\matrixx{0}{1}{-32}{-16}.$$One verifies immediately that $A,
B$, and $C$ belong to $N$, and that their classes in
$N/Z(\BR)\Gamma_0(32)$ have exact orders $2, 4, 4$, respectively.
Thus $T(A), T(B), T(C)$ have exact orders $2, 4, 4$, respectively.
Also $T(B)$ maps $[\infty]$ to itself. Thus $T(B)\in \Aut(E'_\BC)$
is an automorphism of $E'_\BC$ of exact order $4$, proving that
$E'_\BC$ has complex multiplication by $\BZ[i]$. Since $T(A)$ has
fixed point $[i\sqrt{2}/8]$, it is not translation.  It now follows
immediately from that $T(A)^2=1$ that $T(A)=t_{\alpha, -1}$ for some
point $\alpha$ in $E'(\BC)$. Therefore, $T(C)=T(AB^2)=t_{\alpha,
1}$. But $T(C)$ has exact order $4$, whence we see that $\alpha$
must have order $4$. Finally, $T(A)$ is defined over $\BQ$ since it
is the Atkin-Lehner involution. As $T(B^2)$ is the multiplication by
$-1$ and then is clearly defined over $\BQ$. Hence
$\alpha=T(AB^2)([\infty])=[0]$ must be a rational point of exact
order $4$.

Since every elliptic curve over $\BQ$ is known to be parametrized by
the modular curve of the same level as its conductor, it follows
that $E'=(X_0(32), [\infty])$ must be isogenous to the elliptic
curve $y^2=x^3+4x$. However, there are just two isomorphism classes
of curves defined over $\BQ$ in the isogeny class of $E'$, and
$y^2=x^3+4x$ is the unique one with a rational point of order $4$.
Thus $E'$ must be isomorphic to $y^2=x^3+4x$ over $\BQ$. Since
rational points on $y^2=x^3+4x$ of exact order $4$ are $(2, \pm 4)$,
the isomorphism is unique if we require the order 4 point $[0]$ on
$E'$ is mapped to $(2, 4)$. Thus $E'$ has Weierstrass equation
$y^2=x^3+4x$ for unique modular functions $x, y$ such that the cusp
$[0]$ has coordinate $(2, 4)$.

Let $\psi$ denote the unique Gorssencharacter of any elliptic curve
defined over $\BQ$ with conductor 32 and complex multiplication by
$\BZ[i]$. Then the conductor of $\psi$ must be $(1+i)^3\BZ[i]$
because the norm of this conductor times the absolute value $4$ of
the discriminant of $\BQ(i)$ must be 32. It then follows easily from
the main theorem of complex multiplication that
$E'[(1+i)^3]=E'(\BQ(i))_\tor$. But we know that $S\subset
E'(\BQ(i))_\tor$ and has cardinality $8$, thus
$S=E'(\BQ(i))_\tor=E'[(1+i)^3]$.
\end{proof}

For any field extension $F$ over $\BQ$, let $X_0(32)_F$ be the base
change of $X_0(32)$ to $F$ and write $\Aut(X_0(32)_F)$ for the group
of automorphisms of $X_0(32)_F$. Similarly one defines $\Aut(E'_F)$.
Then it is easy to see that
$$\Aut(X_0(32)_F)\cong E'(F) \rtimes \Aut(E'_F).$$
Using the notations in the proof of Proposition \ref{Add 1}, we have
seen above that there is a natural homomorphism
$$T: N \lra \Aut(X_0(32)(\BC), S)\subset \Aut(X_0(32)_\BC),$$
with kernel $Z(\BR)\Gamma_0(32)$.  Now we have that
$\Aut(E'_\BC)\cong \BZ[i]^\times$ and $T(B)\in \Aut(E'_\BC)$ is of
order $4$. One can see that $T(B)$ maps $(x, y)$ to $(-x, iy)$ by
looking at actions of $T(B)$ at $0$ and $B$ at $[\infty]$: at
$[\infty]$, the differential is represented by $dq$ with $q=e^{2\pi
iz}$. It is clear that $B^* dq=dB^* q=i dq$; at $0\in E'$, the
morphism $(x, y)\mapsto (-x, iy)$ brings the Neron differential
$dx/y$ to $i dx/y$.

\begin{prop}\label{Main Prop} With the notations above,
the normalizer $N$ of $\Gamma_0(32)$ is generated by
$Z(\BR)\Gamma_0(32), A$ and $B$.  The homomorphism $T$  induces an
isomorphism
$$N/Z(\BR)\Gamma_0(32) \stackrel{\sim}{\lra} \Aut
(X_0(32)_{\BQ(i)})\cong E'(\BQ(i))\rtimes\Aut(E'_{\BQ(i)}).
$$Moreover, if  write $t_\alpha\in
\Aut(E'_\BC)$ for the translation by $\alpha\in E'(\BC)$, then the
following relations hold:
$$\begin{aligned}
 &t_{(2,
4)}=T\matrixx{0}{1}{-32}{-16}, \quad && t_{(2,
-4)}=T\matrixx{-16}{-1}{32}{0},\quad t_{(0,
0)}=T\matrixx{-2}{-1}{32}{14}, \\
&t_{(-2, 4i)}=T\matrixx{-24}{-7}{32}{8}, \quad && t_{(-2,
-4i)}=T\matrixx{8}{7}{-32}{-24}, \\
& t_{(2i, 0)}=T\matrixx{-4}{-3}{32}{20},\quad && t_{(-2i,
0)}=T\matrixx{4}{1}{32}{12}.
\end{aligned}$$
\end{prop}

\begin{proof}  Since $E'(\BQ)$ has rank 0, we see that $E'(\BQ(i))=E'[(1+i)^3]$
consists of the following $8$ points: $$[\infty], \quad (0, 0),\quad
(2, \pm 4), \quad (\pm 2i, 0), \quad (-2, \pm 4i).$$Note that
$T(C)=t_{(2, 4)}$ and $T(B)$ generate $E'(\BQ(i))\rtimes \Aut
(E'_{\BQ(i)})$. It follows that the image of $T$ contains
$\Aut(X_0(32)_{\BQ(i)})$.  But any $t$ in the image of $T$,
$t([\infty])\in S=E'(\BQ(i))$. It follows that $\Im(T)\subseteq
E'(\BQ(i))\rtimes \Aut(E'_{\BQ(i)})$. Thus the image of $T$ is
$\Aut(X_0(32)_{\BQ(i)}$ and the homomorphism $T$ induces an
isomorphism $N/Z(\BR)\Gamma_0(32) \stackrel{\sim}{\lra} \Aut
(X_0(32)_{\BQ(i)})$. It also follows that $N$ is generated by
$Z(\BR)\Gamma_0(32)$, $A$, and $B$.

Note that $\alpha=(2, 4)$ and $i\alpha$ generate $E'(\BQ(i))$ and
$$t_{-i\alpha}=[-i]\circ [t_\alpha]\circ [i]=T(B^{-1})\circ T(C)\circ
T(B)=T(B^{-1}CB).$$ The verifying of remaining relations  is then
straightforward.
\end{proof}

There is a well known alternative adelic expression for the complex
points of $X_0(32)$, which we will also need.  Let $\BA$ be the
adeles of $\BQ$ and $\BA_f$ its finite part. Let $G=\GL_{2, \BQ}$,
$G(\BA_f)$ its finite-adelic points, and $U_0(32)\subset G(\BA_f)$
the open compact subgroup defined by
$$U_0(32)=\left\{ \begin{pmatrix} a& b\\ c& d\end{pmatrix} \in
\GL_2(\wh{\BZ})\ \Big|\ c\equiv 0\mod 32\wh{\BZ}\right\}.$$ The
complex uniformzation of $X_0(32)$ has the following adelic form
$$X_0(32)(\BC)=G(\BQ)_+\bs (\CH\cup \BP^1(\BQ)) \times G(\BA_f)
/U_0(32).$$For any $z\in \CH \cup \BP^1(\BQ)$ and any $g\in
G(\BA_f)$. we denote by $[z, g]$ its image in $X_0(32)(\BC)$.

The matrix $B^2\in G(\BQ)\subset G(\BA_f)$ normalizes both
$\Gamma=\Gamma_0(32)$ and $U_0(32)$. The morphism $T(B^2)$ is
represented by the Hecke action
$$[z, g]\mapsto [z, gB^{-2}], \quad \forall z\in \CH, g\in G(\BA_f)$$
which is defined over $\BQ$ by the functriality of canonical models
of Shimura varieties. However, the matrix $B$ does not normalize
$U_0(32)$, though it normalizes $\Gamma_0(32)$. The morphism $T(B)$
on $X_0(32)$ can be written as
$$[z, \gamma]\mapsto [z, \gamma B^{-1}], \quad \forall \gamma \in G(\BQ)_+, z\in \CH,$$
but we can not conclude that it is defined over $\BQ$. In fact, it
is defined over $\BQ(i)$.

\bigskip

We now construct suitable points  on $E$ from CM points on
$X_0(32)$. We first consider the case with $n\equiv 1\mod 4$ and the
case with $n\equiv 3\mod 4$ will be considered  later in Theorem
\ref{Main Theorem 2 in section 2}. Note that the set of torsion
points with exact order $4$ on $E$ is the union of the following
subsets:
$$(i, 1-i)+E[2], \quad \left(1+\sqrt{2}, 2+\sqrt{2}\right)+E[2],
\quad\text{and}\quad \left(-1-\sqrt{2},
i\big(2+\sqrt{2}\big)\right)+E[2],$$whose doubles are $(0, 0), (1,
0)$, and $(-1, 0)$, respectively.

\begin{defn} Let $n\equiv 1\mod 4$
be a positive integer and $K=\BQ(\sqrt{-2n})$. Let $P\in
X_0(32)(K^{ab})$ be the image of $
\displaystyle{\frac{i\sqrt{2n}}{8}}$ under the complex
uniformization $\CH \ra X_0(32)$. Define the CM point on $E$
$$z:=f(P)+(1+\sqrt{2}, 2+\sqrt{2})\in E(K^{ab}).$$
For each $t\in \wh{K}^\times$, let $z_t$ denote the Galois
conjugation $z^{\sigma_t}$ of $z$.
\end{defn}

\begin{thm}\label{Main Theorem of section 2} Assume that $n\equiv 1\mod 4$ is a
positive integer. Then, for each  $t\in \wh{K}^\times$, we have
\begin{enumerate}
\item the point $z_t$ is defined over the Hilbert class field $H$ of $K$
and only depends on the class of $t$ modulo $K^\times
\wh{\CO}_K^\times$;
\item the complex conjugation  of $z_t$, denoted by $\bar{z}_t$,  is equal to
$z_{t^{-1}}$; and \item $z_{\varpi t}+z_t=0$ or $(0, 0)$ according
to $n\equiv 1 \mod 8$ or $\equiv 5 \mod 8$.
\end{enumerate}
\end{thm}

\begin{remark} The CM points $z_t$  above are essentially
the same as those Monsky studied in \cite{Monsky} using modular
functions on $X(8)$. Theorem\ \ref{Main Theorem of section 2} still
holds if we replace $z$ by any CM point $\pm f(P)+Q$ where $Q$ is
any 4-torsion point of $E$ with $2Q=(1, 0)$.

\end{remark}
We will show Theorem\ \ref{Main Theorem of section 2} by showing the
following corresponding result on $X_0(32)$ via the modular
parametrization $f$.

\begin{prop} Let $n\equiv 1\mod 4$ be a positive integer and $K=\BQ(\sqrt{-2n})$.
Let $P\in X_0(32)$ be the point defined by $i\sqrt{2n}/8\in \CH$ via
the complex uniformzation.  Let $H'$ be the defining field of $P$.
The the following hold: \begin{enumerate} \item the field $H'\subset
K^{ab}$ of $P$ over $K$ is characterized by
$$\Gal(H'/K)\stackrel{\sim}{\lra}\wh{K}^\times /K^\times (\BZ_2^\times (1+4\CO_2))
\wh{\CO}_K^{\times (2)}$$ via Artin reciprocity law. Here the
supscript in $\wh{\CO}_K^{\times (2)}$ means the component at the
unique place of $K$ above 2 is removed.\item The extension $H'/K$ is
anticyclotomic in the sense that $H'$ is Galois over $\BQ$ such that
the nontrivial involution on $K$ over $\BQ$ acts on $\Gal(H'/K)$ by
the inverse.\item The field $H'$ is a cyclic extension of degree 4
over $H$ with $\Gal(H'/H)$ generated by $\sigma_{1+\varpi}$, where
recall that $\varpi\in K_2^\times$ is the uniformizor $\sqrt{-2n}$.
\item Moreover,
$$P^{\sigma_{1+\varpi}}=P+(-2i, 0), \quad P^{\sigma_{\varpi}}+P=(2,
4).$$\end{enumerate}
\end{prop}

\begin{proof} Recall Shimura's
reciprocity law (for example, see \cite{Lang}). Let $w=i\sqrt{2n}\in
K^\times$ and view $K^\times$ as a sub-torus of $\GL_{2, \BQ}$ via
the $\BQ$-embedding of $K^\times$ into $\GL_{2, \BQ}$: $a+bw\mapsto
\begin{pmatrix} a& -2n b\\ b& a\end{pmatrix}$. Then $w\in \CH$ is the
unique point on $\CH$ fixed by $K^\times$.  For any point
$$x=[w, g]\in X_0(32)(\BC)=G(\BQ)_+\bs (\CH \cup \BP^1(\BQ)) \times
G(\BA_f)/U_0(32), \quad g\in G(\BA_f)$$ and any $t\in
\wh{K}^\times\subset G(\BA_f)$, the action of $\sigma_t$ on $x$ is
given by: $[w, g]^{\sigma_t}=[w, tg]$. It follows that the defining
field $K(x)$ of $x$ is characterized by
$$\Gal(K(x)/K)\simeq \wh{K}^\times/K^\times (\wh{K}^\times \cap
gU_0(32) g^{-1})$$ via the reciprocity law.

Write the point $P=\tiny{\left[w,
\begin{pmatrix} 8&0
\\ 0&1\end{pmatrix}\right]}\in X_0(32)$ in adelic form. Then $H'$
corresponds to the open compact subgroup
$$\wh{K}^\times \cap \matrixx{8}{0}{0}{1} U_0(32)
\matrixx{8^{-1}}{0}{0}{1}=\BZ_2^\times(1+4\CO_2) \wh{\CO}_K^{\times
(2)}.$$It gives the statement (1) and
$$\Gal(H'/H)\stackrel{\sim}{\lra} K^\times \wh{\CO}_K^\times/ K^\times (\BZ_2^\times (1+4\CO_2))
\wh{\CO}_K^{\times (2)}=\CO_2^\times/\BZ_2^\times
(1+4\CO_2)=(1+\varpi)^{\BZ/4\BZ}.$$ Here we use the fact that
$\CO_2^\times=\{\pm 1\}\times 5^{\BZ_2}\times (1+\varpi)^{\BZ_2}$ as
a $\BZ_2$-module. Moreover, since
$$\wh{\BQ}^\times\subset K^\times \cdot(\BZ_2^\times(1+4\CO_2) \wh{\CO}_K^{\times
(2)}),$$the non-trivial involution of $K$ acts on $\Gal(H'/K)$ by
the inverse. The statements (2) and (3) are now proved.

By Proposition \ref{Main Prop}, $t_{(-2i,
0)}=T\tiny{\matrixx{4}{1}{32}{12}}$. Note that $(-2i, 0)$ is of
order $2$.  Thus the relation $P^{\sigma_{1+\varpi}}=P+(-2i, 0)$ is
equivalent to
$$P^{\sigma_{1+\varpi}}=T\matrixx{4}{1}{32}{12}
P,$$ which is just
\begin{eqnarray}\left[w, (1+\varpi)\begin{pmatrix} 8&0\\
0&1\end{pmatrix}\right]=
\left[w, \begin{pmatrix} 8&0\\
0&1\end{pmatrix}\begin{pmatrix} 4&1\\
32&12\end{pmatrix}\right]\end{eqnarray} It is further equivalent to
$$(1+\varpi)\begin{pmatrix} 8&0\\
0&1\end{pmatrix}\in K^\times \begin{pmatrix} 8&0\\
0&1\end{pmatrix}\begin{pmatrix} 4&1\\
32&12\end{pmatrix} U_0(32),$$and then to
$$(1+\varpi)\in K^\times (V\cap \wh{K}^\times), \quad
V=\matrixx{1}{2}{1}{3}\matrixx{8}{0}{0}{1}U_0(32)
\matrixx{8^{-1}}{0}{0}{1}.$$It follows easily from  $n\equiv 1\mod
4$ that
$$V\cap \wh{K}^\times =\wh{\CO}_K^{\times (2)} \BZ_2^\times
(1+\varpi+4\CO_2).$$ Then (2.1) is equivalent to
$$1+\varpi\in K^\times\wh{\CO}_K^{\times (2)} \BZ_2^\times
(1+\varpi+4\CO_2),$$ which is obvious.

By proposition \ref{Main Prop}, $t_{(2,
4)}=T\tiny{\matrixx{-16}{-1}{32}{0}}^{-1}$, thus the relation
$P^{\sigma_\varpi}+P=(2, 4)$ is equivalent to
$$P^{\sigma_\varpi}=T\matrixx{-16}{-1}{32}{0}^{-1} (T(B^2)P),$$
which is just
\begin{eqnarray}
\left[w, \varpi\begin{pmatrix} 8&0\\
0&1\end{pmatrix}\right] &=&
\left[w, \begin{pmatrix} 8&0\\
0&1\end{pmatrix}\begin{pmatrix} 1&1/2\\
0&1\end{pmatrix}\begin{pmatrix} -16&-1\\
32&0\end{pmatrix}\right]
\end{eqnarray}
It is further equivalent to
$$\varpi\matrixx{8}{0}{0}{1}\in K^\times
\matrixx{8}{0}{0}{1}\matrixx{0}{-1}{32}{0}U_0(32)$$and then to
$$\varpi\in K^\times (V\cap \wh{K}^\times), \quad
V=\matrixx{0}{-2}{1}{0} \matrixx{8}{0}{0}{1} U_0(32)
\matrixx{8^{-1}}{0}{0}{1}.$$ It is easy to have that
$$V\cap \wh{K}^\times=\wh{\CO}_K^{\times (2)}\BZ_2^\times
(\varpi+8\CO_2).$$ Then (2.2) is equivalent to $\varpi \in K^\times
\wh{\CO}_K^{\times (2)}\BZ_2^\times (\varpi+8\CO_2)$, which is
obvious. The proof of (4) is complete.
\end{proof}

\begin{proof}[Proof of Theorem \ \ref{Main Theorem of section 2}]
Recall that the defining field $H'=K(P)$ of $P$ over $K$ corresponds
to the subgroup $K^\times \wh{\CO}_K^{\times (2)} \BZ_2^\times
(1+4\CO_2)\subset \wh{K}^\times$. The norm of this subgroup over the
extension $K/\BQ$ is $\BQ^\times \wh{Z}^{\times (2)}\BZ_2^{\times
2}(1+8\BZ_2)$, which corresponds to the abelian extension
$\BQ(\sqrt{2}, i)$ over $\BQ$. Thus $\sqrt{2}, i\in H'$.

For any $t\in \wh{K}^\times$, let
$P_t=\tiny{\left[w, t\begin{pmatrix} 8& \\
&1\end{pmatrix}\right]}=P_1^{\sigma_t}\in X_0(32)$, then
$$z_t:=z^{\sigma_t}=f(P_t)+(1+\sqrt{2}, 2+\sqrt{2})^{\sigma_t}\in E(H').$$
Since $H'/K$ is anti-cyclotomic and $P\in X_0(32)(\BR)$, the complex
conjugation of $f(P_t)$ is $f(P_{t^{-1}})$ and therefore the complex
conjugation of $z_t$ is equal to $z_{t^{-1}}$. This proves (2).

To show (1), we only need to consider the case with $t=1$, i.e.
$z:=z_1\in E(H)$. Note that
$$\sigma_{1+\varpi}(\sqrt{2})=[1+\varpi, K^{ab}/K](\sqrt{2})=[(1+2n)_2,
\BQ^{ab}/\BQ](\sqrt{2})=-\sqrt{2}.$$ Since $\sigma_{1+\varpi}$
generates $\Gal(H'/H)$,  $z\in E(H)$ is equivalent to the relation
$z^{\sigma_{1+\varpi}}=z$, and therefore  is equivalent to
$$f(P^{\sigma_{1+\varpi}})=f(P)+(1+\sqrt{2}, 2+\sqrt{2})-(1+\sqrt{2},
2+\sqrt{2})^{\sigma_{1+\varpi}}=f(P)+(0, 0)=f(P+(-2i, 0))$$ which
follows from the first equality in Proposition 2.6 (4).

To show (3), we only need to show that $z_\varpi+z=0$ or $(0, 0)$
according to $n\equiv 1\mod 8$ or $\equiv 5\mod 8$. Note that
$$[\varpi, K^{ab}/K] (\sqrt{2})=[(2n)_2,
\BQ^{ab}/\BQ](\sqrt{2})=[n_2, \BQ^{ab}/\BQ](\sqrt{2})=(-1)^{(n-1)/4}
\sqrt{2}$$and that $f((2, 4))=(1, 0)$. It follows that
$$\begin{aligned}
z^{\sigma_\varpi}+z&=f(P^{\sigma_\varpi}+P)+\begin{cases}
(1, 0), \quad &\text{if $n\equiv 1\mod 8$},\\
(-1, 0), &\text{if $n\equiv 5\mod 8$}.
\end{cases}\\
&=f(P^{\sigma_\varpi}+P-(2, 4))+\begin{cases}
0, \quad &\text{if $n\equiv 1\mod 8$},\\
(0, 0), &\text{if $n\equiv 5\mod 8$}.
\end{cases}
\end{aligned}$$Thus the desired follows then from the second equality in
Proposition 2.6 (4).

\end{proof}

We now consider the case with $n\equiv 3\mod 4$.
\begin{defn}
Let  $n\equiv 3\mod 4$ be a positive integer and
$K=\BQ(\sqrt{-2n})$. Let $P\in X_0(32)$ be the image of
$\displaystyle{\frac{2+i\sqrt{2n}}{8}}$ under the complex
uniformization and define a CM point
$$z:=f(P)+(1+\sqrt{2}, 2+\sqrt{2})\in E(K^{ab}).$$
For any $t\in \wh{K}^\times$ satisfying $\sigma_t$ fixing $i$,
define $z_t:=z^{\sigma_t}$.
\end{defn}

\begin{thm}  \label{Main Theorem 2 in section 2} Assume that $n\equiv 3\mod
4$ is a positive integer. Then, for each $t\in \wh{K}^\times$ with
$\sigma_t$ fixing $i$, we have
\begin{enumerate}
\item the point $z_t\in E(H(i))$ and the involution $\sigma_{1+\varpi}$ of $\Gal(H(i)/H)$ maps $z_t$ to
$z_t+(0, 0)$;
\item the complex conjugation
 of $z_t$, denoted by $\bar{z}_t$, is equal to $-z_{t^{-1}}+(1, 0)$;
\item let $\varpi'=\varpi(1+\varpi)\in K_2^\times$ (so that
$\sigma_{\varpi'}$ fixes $i$), then $z_{\varpi' t}-z_t=(1, 0)$ or
$(-1, 0)$ according to $n\equiv 7\mod 8$ or $3\mod 8$.
\end{enumerate}
\end{thm}
We will give the proof of Theorem \ref{Main Theorem 2 in section 2}
after we prove the following
\begin{prop} Let $n\equiv 3\mod 4$ be a positive integer and $P\in
X_0(32)$ be the CM point corresponding
$\displaystyle{\frac{2+i\sqrt{2n}}{8}\in \CH}$ via complex
uniformzation. The defining field $H'\subset K^{ab}$ of $P$ over $K$
is characterized by
$$\Gal(H'/K)\stackrel{\sim}{\lra}\wh{K}^\times /K^\times (\BZ_2^\times (1+4\CO_2))
\wh{\CO}_K^{\times (2)}$$ via Artin reciprocity law so that
$\Gal(H'/H)$ is generated by $\sigma_{1+\varpi}$. Moreover,
$$P^{\sigma_{1+\varpi}}=P+(2i, 0), \qquad P^{\sigma_{\varpi'}}-P=(2,
4).$$Here $\varpi'=\varpi(1+\varpi)\in K_2^\times$ so that
$\sigma_{\varpi'}$ fixes $i$.
\end{prop}

\begin{proof}The proof of the first part is the same as in the case
$n\equiv 1\mod 4$. By Proposition 2.2, $t_{(2i,
0)}=T\matrixx{-4}{-3}{32}{20}$. The relation
$P^{\sigma_{1+\varpi}}=P+(2i, 0)$ is equivalent to
$$P^{\sigma_{1+\varpi}}=T\matrixx{-4}{-3}{32}{20}P$$
which is
\begin{eqnarray}\left[w,
(1+\varpi)\matrixx{8}{-2}{0}{1}\right]=\left[w,
\matrixx{8}{-2}{0}{1}\matrixx{-4}{-3}{32}{20}^{-1}\right].
\end{eqnarray} It is further equivalent to,
$$(1+\varpi)\matrixx{8}{-2}{0}{1}\in K^\times
\matrixx{8}{-2}{0}{1}\matrixx{20}{3}{-32}{-4} U_0(32),$$ and then to
$$(1+\varpi)\in K^\times (V\cap \wh{K}^\times), \quad V=
\matrixx{7}{8}{-1}{-1}\matrixx{8}{0}{0}{1}U_0(32)\matrixx{8}{0}{0}{1}^{-1}\matrixx{1}{-2}{0}{1}^{-1}.$$
But it is easy to see that $V\cap \wh{K}^\times=\wh{\CO}_K^{\times
(2)} \BZ_2^\times (1+\varpi+4\CO_2)$ provided that $n\equiv 3\mod
4$. Therefore (2.3) is equivalent to
$$1+\varpi\in K^\times\wh{\CO}_K^{\times (2)} \BZ_2^\times
(1+\varpi+4\CO_2),$$ which is obvious.

By Proposition 2.2, $t_{(2, 4)}=T\matrixx{0}{1}{-32}{-16}$. Thus the
relation $P^{\sigma_{\varpi'}}-P=(2, 4)$ is equivalent to
\begin{eqnarray}\left[w,
\varpi(1+\varpi)\matrixx{8}{-2}{0}{1}\right]=\left[w,
\matrixx{8}{-2}{0}{1}\matrixx{0}{1}{-32}{-16}^{-1}\right].
\end{eqnarray} It is further equivalent to,
$$\varpi(1+\varpi)\matrixx{8}{-2}{0}{1}\in K^\times
\matrixx{8}{-2}{0}{1}\matrixx{-16}{-1}{32}{0} U_0(32),$$ and then to
$$\varpi(1+\varpi)\in K^\times (V\cap \wh{K}^\times), \quad V=
\matrixx{-6}{-2}{1}{0}\matrixx{8}{0}{0}{1}U_0(32)\matrixx{8}{0}{0}{1}^{-1}\matrixx{1}{-2}{0}{1}^{-1}.$$
But it is easy to see that $V\cap \wh{K}^\times=\wh{\CO}_K^{\times
(2)} \BZ_2^\times (2+\varpi+8\CO_2)$. Therefore (2.4) is equivalent
to
$$\varpi(1+\varpi)\in K^\times\wh{\CO}_K^{\times (2)} \BZ_2^\times
(2+\varpi+8\CO_2),$$ which is obvious since $n\equiv 3\mod 4$.

\end{proof}

\begin{proof}[Proof of Theorem \ref{Main Theorem 2 in section 2}]
It is clear that $\sqrt{2}\in H$ but $i\notin H$ and that the Galois
group $\Gal(H'/H(i))=\{1, \sigma_{1+\varpi}^2\}$. Thus the item (1)
follows from the relation $z^{\sigma_{1+\varpi}}=z+(0, 0)$. It is
clearly that the relation is equivalent to
$f(P^{\sigma_{1+\varpi}}-P-(2i, 0))=0$ by noting $f((2i, 0))=(0,
0)$, which is given in Proposition 2.9.

Note that the  point $[i\sqrt{2n}/8]$ is real. Then the complex
conjugation  of $f(P)$ is $-f(P)$ and therefore the complex
conjugation $\ov{z}$ of $z$ is equal to $-z+(1, 0)$. Thus we have
that for any $t\in \wh{K}^\times$ fixing $i$,
$\ov{z_t}=\ov{z}^{t^{-1}}=-z_{t^{-1}}+(1, 0)$. This proves the item
(2).

For item (3), it is enough to show the case with $t=1$. Now
$\sigma_{\varpi'}(\sqrt{2})=\sigma_\varpi(\sqrt{2})=[(2n)_2,
\BQ^{ab}/\BQ](\sqrt{2})=(-1)^{(n^2-1)/8}\sqrt{2}$.  By the relation
$P^{\sigma_{\varpi'}}-P=(2, 4)$ and $f((2, 4))=(1, 0)$, we have
$$\begin{aligned}
z_{\varpi'}-z&=f(P^{\sigma_{\varpi'}}-P)+(1+\sqrt{2},
2+\sqrt{2})^{\sigma_\varpi}-(1+\sqrt{2}, 2+\sqrt{2})\\
&=(1, 0)+\begin{cases} (0, 0), \quad &\text{if $n\equiv 3\mod
8$},\\
0, &\text{if $n\equiv 7\mod 8$},
\end{cases}
\end{aligned}$$which is equal to $(-1,0)$ or $(1, 0)$ according to
$n\equiv 3\mod 8$ or $\equiv 7\mod 8$,
\end{proof}

\section{Comparsion of Heegner Points}
Let $n=p_0p_1\cdots p_k$ be a product of distinct odd primes with
$p_i\equiv 1\mod 8$ for $1\leq i\leq k$ and $p_0\nequiv 1\mod 8$.
Let $m_0$ be a positive divisor of $2n$ such that $m_0\equiv 5, 6$,
or $7\mod 8$. In this section, we will generalize the construction
in Theorem \ref{Main Theorem 3} to define a point $P^\chi(f)\in
E(\BQ(\sqrt{m_0^*}))^-$, where $m_0^*=(-1)^{(n-1)/2} m_0$. With $n$
replaced by the odd part $n_0$ of $m_0$, this construction gives a
point $P^{\chi_0}(f)$ already define in Theorem \ref{Main Theorem 3}
for $n=n_0$. By Gross-Zagier and Kolyvagin, these two points are
actually linearly dependent modulo $E[2]$. The main result of this
section is a comparison of 2-divisibility of these two points using
the generalized Gross-Zagier formula and the 2-divisibility of
special values of L-series. The comparison will be a key ingredient
in the induction argument of proving Theorem \ref{Main Theorem 3} in
the next section. Let us begin with some notations.

Put
$$K=\BQ(\sqrt{-2n}), \qquad K_0=\BQ(\sqrt{-2n_0}),$$
and write $\eta, \eta_0$ for the abelian characters of $\BQ$
defining these two quadratic fields.  Let
$m_0^*=(-1)^{\frac{n_0-1}{2}}m_0$  and we will also consider the two
 extensions
$$J=K(\sqrt{m_0^*}), \qquad J_0=K_0(\sqrt{m_0^*})$$
and write $\chi, \chi_0$ for the abelian characters of $K$ and
$K_0$, respectively, defining them. These two extensions are both
unramified. In fact, it is easy to see that they are contained in
genus subfields of $K, K_0$, respectively.  For any non-zero integer
$d$, let us write $L(E^{(d)}, s)$ for the complex L-function of the
elliptic curve $E^{(d)}: dy^2=x^3-x$ over $\BQ$. On the other hand,
we write $L(E/K, \chi, s), L(E/K_0, \chi_0, s)$ for the complex
L-functions of $E$ over $K, K_0$, twisted by the characters $\chi,
\chi_0$, respectively.
\begin{lem}Let $c\in \{1, 2\}$ denote the integer $2n_0/m_0$.
Then the following equalities hold:
$$L(E/K, \chi, s)=L(E^{(cn/n_0)}, s)L(E^{(m_0)}, s),\quad
L(E/K_0, \chi_0, s)=L(E^{(c)}, s)L(E^{(m_0)}, s).$$
\end{lem}
The proof is an immediate consequence of the Artin formalism applied
to the L-functions of $E$ for the extensions $J, J_0$ of $\BQ$,
which are quartic when $m_0^*\neq -2n$. For example, the first
equality follows on noting that the induced character of $\chi$ is
the sum of the characters defining the two quadratic extensions
$\BQ(\sqrt{m_0^*})$ and $\BQ(\sqrt{-2n/m_0^*})$. Note also that, for
any non-zero integer $d$, the curves $E^{(d)}$ and $E^{(-d)}$ are
isomorphic over $\BQ$. We hope that the usefulness of such a Lemma
in an inductive argument is immediately clear. Indeed, it is very
well known that $L(E^{(c)}, s)$ does not vanish at $s=1$. Thus
$L(E/K, \chi, s)$ and $L(E/K_0, \chi_0, s)$ will have a zero of the
same order at $s=1$ if and only if $L(E^{(cn/n_0)}, s)$ does not
vanish at $s=1$.  We note that this latter assertion does not always
hold. For example, if we take $n=p_0p_1$, $n_0=p_0\equiv 5\mod 8$ to
be any prime,  and $p_1=17$, then $2n/n_0=34$ and $L(E^{(34)}, s)$
has a zero of order 2 at $s=1$ (indeed, 34 is the smallest square
free congruent number which does not lie in the residue classes of
$5, 6, 7\mod 8$).  Nevertheless, the work of Zhao (See Proposition
\ref{Prop C}) always provides a lower bound of the power of $2$
dividing the algebraic part of  $L(E^{(cn/n_0)}, 1)$, which is
precisely what we will need to carry out an induction argument on
$k$, via the comparison of the heights of two Heegner points on $E$,
 which we now construct.

View $K$ as a $\BQ$-subalgebra of $M_2(\BQ)$ via the embedding
$$K\lra M_{2\times 2}(\BQ), \quad a+b\sqrt{-2n}/8 \lto
\matrixx{a}{-2nb/64}{b}{a}, \quad \forall a, b\in \BQ,$$with which
$K^\times$ is a $\BQ$-subtorus of $\GL_{2, \BQ}$ such that
$i\sqrt{2n}/8$ is the unique fixed point of $K^\times$ on the upper
half complex plane $\CH$. Define a CM point $P\in X_0(32)$ to be
$[w/8, 1]$ if $n\equiv 1\mod 4$ and $[(w+2)/8, 1]$ if $n\equiv 3\mod
4$. Let $f: X_0(32)\ra E$ be a modular parametrization of degree 2
of the elliptic curve $E: y^2=x^3-x$.  By Theorem \ref{Main Theorem
of section 2} and \ref{Main Theorem 2 in section 2}, we have that
$f(P)+(1+\sqrt{2}, 2+\sqrt{2})\in E(H(i))$. Thus $f(P)\in E(H(i))$
since $\sqrt{2}\in H(i)$.  Here $H$ is the Hilbert class field of
$K$. Define
$$P^\chi(f)=\sum_{\sigma\in \Gal(H(i)/K)} f(P)^\sigma \chi(\sigma)
\in E(K(\sqrt{m_0^*}))^-,$$where $E(K(\sqrt{m_0^*}))^-$ is the
subgroup of points in $E(K(\sqrt{m_0^*}))$ which are mapped to their
inverses under the involution of $K(\sqrt{m_0^*})$ over $K$.
Similarly, let $E(\BQ(\sqrt{m_0^*}))^-$ denote the subgroup of
points in $E(\BQ(\sqrt{m_0^*}))$ which are mapped to their inverses
under the non-trivial involution of $\BQ(\sqrt{m_0^*})$ over $\BQ$.
Note that  $\chi$ in this section is the character defining the
extension $K(\sqrt{m_0^*})$ over $K$, but is not the one in the
introduction defining $K(\sqrt{m^*})$ if $m_0\neq m$.
\begin{lem} The point  $P^\chi(f)$  belongs to
$E(\BQ(\sqrt{m_0^*}))^-$.
\end{lem}
\begin{proof}  Note that when $m_0^*\neq -2n$, the extension $K(\sqrt{m_0^*})$
over $\BQ$ is quartic and its Galois group is generated by complex
conjugation and the non-trivial element in
$\Gal(K(\sqrt{m_0^*})/K)$. In this case,  we only need to check that
$P^\chi(f)\in E(\BQ(\sqrt{m_0^*}))$.

If $n\equiv 5\mod 8$, then $m_0^*=m_0=n_0$ is positive. Note that
$P=[i\sqrt{2n}/8, 1]$ is defined over $\BR$ since $i\sqrt{2n}$ is
pure imaginary.  Note also that  the complex conjugation acts on
$\Gal(H(i)/K)$ by the inverse. It follows that $P^\chi(f)$ is
invariant under complex conjugation, and therefore belongs to
$E(\BQ(\sqrt{m_0^*}))$.

Now assume $n\equiv 3\mod 4$, then $m_0^*=-m_0$.  Note that $P$ is
the multiplication of a real point by $[i]$ and therefore is mapped
to its negative under the complex conjugation. It follows that
$P^\chi(f)$ is mapped to $-P^\chi(f)$ under the complex conjugation.
If $m_0=2n$, it is now clear that $P^\chi(f)\in
E(\BQ(\sqrt{-m_0}))^-$. When $m_0\neq 2n$, choose $\sigma \in
\Gal(H(i)/K(i))$ mapping $\sqrt{-m_0}$ to $-\sqrt{-m_0}$, then both
$\sigma$ and the complex conjugation take $P^\chi(f)$ to
$-P^\chi(f)$. Therefore their composition fixes $P^\chi(f)$ and has
fixed field $\BQ(\sqrt{-m_0})$ in $K(\sqrt{-m_0})$. This shows that
$P^\chi(f)\in E(\BQ(\sqrt{-m_0}))$.
\end{proof}

Analogously for $K_0$, let $P_0\in X_0(32)(K_0^{ab})$ be the point
$[\sqrt{-2n_0}/8, 1]$ if $n_0\equiv 1\mod 4$ and
$[(\sqrt{-2n_0}+2)/8, 1]$ if $n_0\equiv 3\mod 4$. Replacing $K, P,
\chi$ by $K_0, P_0, \chi_0$, we similarly obtain another point
$P^{\chi_0}(f)\in E(\BQ(\sqrt{m_0^*}))^-$.

The main goal of this section is to compare  $2$-divisibilities of
the two points $P^\chi(f)$ and $P^{\chi_0}(f)$ in the group
$E(\BQ(\sqrt{m_0^*}))^-$, which is given by the following result.
\begin{thm}\label{Main Theorem in section 3}  Let $n=p_0p_1\cdots p_k$ and $n_0$,  $m$, and $m_0$
be integers as above. Assume that $n_0=p_0p_{i_1}\cdots p_{i_{k-s}}$
is a proper divisor of $n$, i.e. $s>0$, and that $P^{\chi_0}(f)$
belongs to $2^{t} E(\BQ(\sqrt{m_0^*}))^-+E[2]$ for some integer
$t\geq 0$. Then
$$P^\chi(f)\in 2^{t+s+1}E(\BQ(\sqrt{m_0^*}))^-+E[2].$$
\end{thm}

The proof of Theorem \ref{Main Theorem in section 3} is divided into
three steps. First, reduce the comparison of the two points in
Mordell-Weil group to the comparison of their heights via
Kolyvagin's result; second, further reduce to the comparison of two
Special L-values via generalized Gross-Zagier formula; third,
estimate 2-adic valuations of these special L-values.

\begin{prop}\label{Prop A} If either $P^{\chi_0}(f)$ or $P^\chi(f)$ is not torsion, then
$E(\BQ(\sqrt{m_0^*}))^-\otimes_\BZ \BQ$ is  one dimensional
$\BQ$-vector space. In this case,  the ratio $[P^\chi(f):
P^{\chi_0}(f)]\in \BQ\cup \{\infty\}$ of the two points  in this one
dimensional space is given by

\begin{eqnarray}
[P^\chi(f): P^{\chi_0}(f)]^2=[\wh{h} (P^\chi(f)):
\wh{h}(P^{\chi_0}(f))],
\end{eqnarray} where $\wh{h}:
E(\bar{\BQ})\ra \BR$ denotes the N\'{e}ron-Tate height function.
\end{prop}
\begin{proof}Note that Heegner hypotheis is not satisfied for $(E, K, \chi)$, but
one can still use Kolyvagin's Euler system method to see that
$E(K(\sqrt{m_0^*}))^-$ is of rank one if $P^\chi(f)$ is not torsion
(for example, see Theorem 3.2 of \cite{Nek}) and therefore
$E(\BQ(\sqrt{m_0^*}))^-$ is of rank one by Lemma 3.2. Another
argument with Kolyvagin's original result is as follows. Using the
generalized Gross-Zagier formula (Theorem 1.2 in \cite{YZZ}), we
know the L-series $L(s, E, \chi)=L(s, E^{(2n/m_0)})L(s, E^{(m_0)})$
has vanishing order 1 at the central point $s=1$. Considering
$\epsilon$-factor, we then know that $L(s, E^{(m_0)})$ has vanishing
order $1$ at $s=1$. Taking an imaginary quadratic field $K'$ such
that the Heegner hypothesis is satisfied for $(E^{(m_0)}, K')$ and
$L(s, E^{(m_0)}_{K'})$ has vanishing order 1 at $s=1$, then
Kolyvagin's original result shows that $E^{(m_0)}(\BQ)\cong
E(\BQ(\sqrt{m_0^*}))^-$ is of rank one. Similarly, if
$P^{\chi_0}(f)$ is not torsion then $E(\BQ(\sqrt{m_0^*}))^-$ is of
rank one.

The equality (3.1) now follows from that the height function
$\wh{h}$ is quadratic.
\end{proof}

Using the generalized Gross-Zagier formula (Theorem 1.2 in
\cite{YZZ}), we  further express this ratio in term of special
L-values. For any square-free integer $d\geq 1$, let  $\Omega^{(d)}$
denote the real period of $E^{(d)}$, defined by
$$\Omega^{(d)}:=\frac{2}{\sqrt{d}}\int_1^\infty
\frac{dx}{\sqrt{x^3-x}}.$$It is known that the algebraic part of
L-values
$$L^\alg (E^{(d)}, 1):=L(E^{(d)}, 1)/\Omega^{(d)}$$
is a rational number.  Note that $L^\alg(E^{(c)}, 1)\neq 0$ for
$c=1$ and $2$ (see \cite{BSD} p. 87 ).
\begin{prop}\label{Prop B} We have that
\begin{eqnarray}
\wh{h} (P^\chi(f))=\frac{L^\alg(E^{(cn/n_0)}, 1)} {L^\alg( E^{(c)},
1)}\cdot \wh{h} (P_0^{\chi_0}(f)).
\end{eqnarray}
\end{prop}

\begin{remark}In the proof of this proposition, we will use
the language of automorphic representation.  Let $\pi^{(d)}$ denote
the automorphic representation associated to the elliptic curve
$E^{(d)}$.  Let $L(s, \pi^{(d)})$ denote its complete L-series and
$L^{(\infty)}(s, \pi^{(d)})$ its finite part. Then
$$L^{(\infty)}\big(s, \pi^{(d)}\big)=L\Big(E^{(d)}, s+\frac{1}{2}\Big).$$
Moreover, we have that
$$\frac{L^\alg(E^{(cn/n_0)}, 1)} {L^\alg( E^{(c)},
1)}=\frac{L(1/2, \pi^{(cn/n_0)})/\Omega^{(cn/n_0)}}{L(1/2,
\pi^{(c)})/\Omega^{(c)}}.$$
\end{remark}

Before the proof of this proposition, we need to recall the
generalized Gross-Zagier formula in \cite{YZZ}.

Let $G=\GL_{2, \BQ}$. Let $X=\varprojlim_U X_U$ be the projective
limit of modular curves indexed by open compact subgroup $U\subset
G(\BA_f)$. Let $\xi=(\xi_U)$ be the compatible system of Hodge class
such that $\xi_U$ is represented by $\infty$ on each geometric
irreducible component of $X_U$. For the elliptic curve $E:
y^2=x^3-x$, define
$$\pi_E:=\Hom^0_\xi(X, E):=\varinjlim_U \Hom^0_{\xi_U} (X_U, E)$$
where $\Hom^0_{\xi_U} (X_U, E)$ is the group of the morphisms $f$ in
$\Hom_\BQ(X_U, E)\otimes \BQ$ satisfying $f(\infty)=0$. The
$\BQ$-vector space $\pi_E$ is endowed with  the natural
$G(\BA_f)$-structure. Let $\BH$ be the division quaternion algebra
over $\BR$ and let $\pi_\infty=\BQ$ be the trivial representation of
$\BH^\times$. Then the representation $\pi=\pi_E$ is a restricted
tensor product $\pi=\otimes_{v\leq \infty} \pi_v$ (with respect to a
spherical family $f_v^\circ$) as a representation of the incoherent
group $\BG=G(\BA_f)\times \BH^\times$.  We call $\pi_E$ the rational
automoorphic representation associated to the elliptic curve $E$.

The representation $\pi=\pi_E$  is self-dual via the  perfect
$\BG$-invariant pairing
$$(\ ,\ ): \pi\times \pi \lra \BQ$$
given by
$$(f_1, f_2)=\Vol(X_U)^{-1}  f_{1, U}\circ f_{2, U}^\vee,\qquad
\Vol(X_U)=\int_{X_U(\BC)} \frac{dxdy}{2\pi y^2}$$ where $f_i$ is
represented by $f_{i, U}\in \Hom_{\xi_U}(X_U, E)=\Hom^0(J_U, E)$
with $J_U$ the Jacobin of $X_U$, and $f_{2, U}^\vee: E\ra J_U$ is
the dual of $f_{2, U}$ so that $f_{1, U}\circ f_{2, U}^\vee\in \BQ=
\End^0(E)$. For each place $v$, let $(\ ,\ )_v$ be a
$G(\BQ_v)$-invariant pairing such that for almost all unramified
$v\nmid \infty$, $(f_v^0, f_v^0 )_v=1$ and such that for
$f=\otimes_v f_v\in \pi_E$, $(f, f)=\prod_v (f_v , f_v)_v$.

Let $\eta: \wh{\BQ}^\times/\BQ^\times\ra \{\pm 1\}$ be the character
associated to the quadratic extension $K/\BQ$. Let
$K^1=K^\times/\BQ^\times$. Fix a Haar measure $dt_v$ on
$K^1(\BQ_v)=K_v^\times/\BQ_v^\times$ for each place $v$ of $\BQ$
such that the product measure over all $v$ is the Tamagawa measure
on $K^1\bs K^1(\BA)$  multiple by $L(1, \eta)$. Define
$$\beta_v(f_v)=\frac{L(1, \eta_v) L(1, \pi_v, \ad)}{\zeta_{\BQ_v}(2)
L(1/2, \pi_v, \chi_v)}\int_{K_v^\times/\BQ_v^\times} (\pi_v(t)f_v,
f_v)_v \chi_v(t) dt_v.$$

The Gross-Zagier formula (\cite{YZZ} Theorem 1.2)  says:
\begin{eqnarray}\qquad\qquad (2h_K)^{-2}\wh{h}
(P_\chi(f))=\frac{\zeta_\BQ(2) L'(1/2, \pi, \chi)} {4L(1, \eta)^2
L(1, \pi, \ad)} \prod_v \beta_v(f_v),\quad \forall \ f=\otimes_v f_v
\end{eqnarray} where all L-functions, including $\zeta_\BQ$, are all
complete L-series and $L(s, \pi, \chi)$  (resp $L(s, \pi, \ad)$) is
L-series defined for Jacquet-Langlands correspondence of $\pi$, and
$h_K$ is the ideal class number of $K$. Similarly,
\begin{eqnarray}\qquad\qquad (2h_{K_0})^{-2}\wh{h} (P^{\chi_0}(f))=\frac{\zeta_\BQ(2)
L'(1/2, \pi, \chi_0)}{4L(1, \eta_0)^2 L(1, \pi, \ad)} \prod_v
\beta_{0,v}(f_v),\quad \forall \ f=\otimes_v f_v\end{eqnarray} where
the subscript $0$ correspondes to $K_0=\BQ(\sqrt{-2n_0})$.

We use the formulaes (3.3) and (3.4) to compute the ratio of
$P^\chi(f)$ and $P_0^{\chi_0}(f)$. It is more convenient to fix Haar
measures as follows (which are used in \cite{YZZ}).
\begin{itemize}
\item for each place $v$ of $\BQ$, the Haar measure $dx_v$ on $\BQ_v$
is  self dual with respect to the standard additive character
$\psi_v$ on $\BQ_v$: $\psi_\infty(x)=e^{2\pi ix}$ and
$\psi_v(x)=e^{-2\pi i \iota_v(x)}$ for $v\nmid \infty$ where
$\iota_v: \BQ_v/\BZ_v \ra \BQ/\BZ$ is the natural embedding. The
Haar measure $dx_v^\times$ on $\BQ_v^\times$ is given by
$\zeta_{\BQ_v}(1)|x|_v^{-1} dx_v$,  where $|\ |_v$ is the normalized
absolute value on $\BQ_v$.
\item the Haar
measure $dy$ on $K_v$ is such that the Fourier transformation
$$\wh{\Phi}(x)=\int_{K_v} \Phi(y) \psi_v(\pair{x, y}) dy$$
satisfies $\wh{\wh{\Phi}}(x)=\Phi(-x)$ where the pairing is
$$\pair{x,
y}=\RN_{K_v/\BQ_v}(x+y)-\RN_{K_v/\BQ_v}(x)-\RN_{K_v/\BQ_v}(y).$$ The
Haar measure $d^\times y$ on $K_v^\times$ is given by
$$d^\times y=\zeta_{K_v}(1)|\RN_{K_v/\BQ_v}(y)|_v^{-1} dy$$
where $\zeta_{K_v}(s)=\zeta_{\BQ_v}(s)^2$ for $v$ splits in $K$. For
any $v\nmid \infty$, let $D$ be the discriminant of $K_v$ in
$\BZ_p$, then
$$\Vol(\CO_{K_v}, dy)=\Vol(\CO_{K_v^\times}, d^\times
x)=|D|_v^{1/2}.$$
\item take the quotient Haar measure $dt_v$ on
$K_v^1=K_v^\times/\BQ_v^\times$ and then the product Haar measure
$\otimes_v dt_v$ has the total volume of $K^1(\BQ)\bs K^1(\BA)$
equal to $2L(1, \eta)$, which  satisfies the requirement in the
Gross-Zagier formula above. Then
$$\Vol(K_v^1, dt_v)=\begin{cases}
2, \qquad &\text{if $v=\infty$,}\\
1, &\text{if $v$ is inert in $K$,}\\
2|D|_v^{1/2}, &\text{if $K_v/\BQ_p$ is ramified}.
\end{cases}$$
\end{itemize}
\begin{lem}Let $f: X_0(32)\ra E$ be a degree $2$ modular parametrization.
With the above fixed Haar measure,  we have
$$\beta_v(f_v)\big/\beta_{0, v}(f_v)=\begin{cases}
p^{-1/2}, \qquad &\text{if $v=p\Big|\displaystyle{\frac{n}{n_0}}$},\\
1, &\text{otherwise.}
\end{cases}$$
\end{lem}
\begin{proof}
Note that $\beta_v(f_v)/(f_v, f_v)=1$ in the following spherical
case: $K_v/\BQ_v, \pi_v, \chi_v$ are all unramified and $f\in
\pi^{G(\BZ_v)}$ and $\CO_{K_v}^\times/\BZ_v^\times$ has volume one.
Thus
$$\beta_v(f_v)/(f_v, f_v)=1,\ \forall v\nmid 2n \infty, \qquad
\beta_{0, v}(f_v)/(f_v, f_v)=1, \ \forall v\nmid 2n_0 \infty$$ and
therefore $\beta_v(f_v)\big/\beta_{0, v}(f_v)=1$ for any $v\nmid 2n
\infty$.

Now let $p$ be an odd prime. Then $\pi_p$ is unramified and $f_p$ is
a non-zero (spherical) vector in the one-dimensional space
$\pi_p^{G(\BZ_p)}$. Then the normalized matrix coefficient
$$\Psi_p(g):=(\pi_p(g) f_p, f_p)/(f_p, f_p), \qquad g\in G(\BQ_p)$$
 is bi-$G(\BZ_p)$-invariant and satisfies the Macdonald formula (See \cite{Bump} Theorem 4.6.6):
$$\Psi_p\left(\begin{pmatrix} p^m & \\ & 1\end{pmatrix}\right)=\frac{p^{-m/2}}{1+p^{-1}} \left(\alpha^m \frac{1-p^{-1}\alpha^{-2}}{1-\alpha^{-2}}+\alpha^{-m}  \frac{1-p^{-1}\alpha^2}{1-\alpha^2}\right), \quad m\geq 0.$$
Here $(\alpha, \alpha^{-1})$ are the Satake parameteers of $\pi_p$.

For $p|n$ which is then ramified in $K$ and let $\varpi_p$ be a
uniformizor of $K_p$,   using the above formula and the
decomposition
$$K_p^\times/\BQ_p^\times=(\CO_{K_p}^\times/\BZ_p^\times) \cup (\varpi_p\CO_{K_p}^\times/\BZ_p^\times), \qquad \varpi_p=\sqrt{-2n}\in G(\BZ_p) \begin{pmatrix} p&\\ &1\end{pmatrix} G(\BZ_p),$$
we compute the integral
$$\begin{aligned}
\int_{K_p^\times/\BQ_p^\times} \Psi_p(t) \chi_p(t)dt&=\left(1+\Psi_p(\varpi_p)\chi_p(\varpi_p)\right) \Vol(\CO_{K_p}^\times/\BZ_p^\times)\\
&=\frac{\Vol(\CO_{K_p}^\times/\BZ_p^\times)}{1+p^{-1}} \cdot
\left(1+\alpha \chi_p(\varpi_p)p^{-1/2}\right)\left(1+\alpha^{-1}
\chi_p(\varpi_p)p^{-1/2}\right)
\end{aligned}$$
By the following formula for local factors:
$$\begin{aligned}
L(1, \pi_p, \ad)&=(1-\alpha^2 p^{-1})^{-1}(1-\alpha^{-2}p^{-1})^{-1}
(1-p^{-1})^{-1}\\
L(1/2, \pi_p, \chi_p)&=(1-\alpha
\chi_p(\varpi)p^{-1/2})^{-1}(1-\alpha^{-1}\chi_p(\varpi)p^{-1/2})^{-1},
\end{aligned}$$
we have
$$
\frac{\beta_p(f_p)}{(f_p, f_p)}= \frac{L(1, \eta_p) L(1,\pi_p,
\ad)}{\zeta_{\BQ_p}(2)L(1/2, \pi_p, \chi_{p})}\cdot
\int_{K_p^\times/\BQ_p^\times} \Psi_p(t) \chi_p(t)
dt=\Vol(\CO_{K_p}^\times/\BZ_p^\times)=p^{-1/2}.$$

It follows that $\beta_p(f_p)/\beta_{0, p}(f_p)=p^{1/2}$ for each
$p|(n/n_0)$ and $=1$ for all $p\nmid 2\infty n/n_0$. But
$n/n_0\equiv 1\mod 8$ implies that $n/n_0=\gamma^2$ for some
$\gamma\in \BQ_2^{\times}$, i.e. $K_2\simeq K_{0, 2}$. Note that
$f_2$ is a newvector in $\pi_2$, i.e. invariant under the
$2$-component $U_0(32)_2$ of $U_0(32)$ and that the two embeddings
of $K_2^\times=K_{0, 2}^\times$ into $\GL_2(\BQ_2)$ are conjugate by
$\begin{pmatrix}1 &\\ & \gamma
\end{pmatrix}\in U_0(32)_2$. Now it is easy to see that
$\beta_2(f_2)/\beta_{0, 2}(f_2)=1$.  Note that $\pi_\infty=\BQ$ is
the trivial representation of $\BH^\times$, one can also easily
check that $\beta_{\infty}(f_\infty)/\beta_{0, \infty}(f_\infty)=1$.

\end{proof}

\begin{proof}[Proof of Proposition \ref{Prop B}]

Recall that for a non-zero integer square-free $d$, we denote by
$\pi^{(d)}$ the automorphic representation corresponding to the
elliptic curve $E^{(d)}: dy^2=x^3-x$. Then $\pi^{(d)}=\pi^{(-d)}$
and by Lemma 3.1,
$$L(s, \pi, \chi)=L(s, \pi^{(m_0)})L(s, \pi^{(cn/n_0)}),
\quad L(s, \pi, \chi_0)=L(s, \pi^{(m_0)})L(s, \pi^{(c)}).$$ Note
that the functional equation of $L(s, \pi^{(m_0)})$ (resp. $L(s,
\pi^{(2n/m_0)})$)  has sign -1 (resp.  1) since $m_0\equiv m\equiv
5, 6$, or $7\mod 8$ by our assumption.  Thus
$$L'(1/2, \pi, \chi)/L'(1/2, \pi, \chi_0)=L(1/2, \pi^{(cn/n_0)})/L(1/2, \pi^{(c)}).$$

Since the special value $L(1/2, \pi^{(c)})\neq 0$ for $c=1, 2$, thus
 $P_0^{\chi_0}(f)$ is torsion if and only if $L'(1/2,
\pi^{(m_0)})=0$. It follows that if $P^{\chi_0}(f)$ is torsion then
$P^\chi(f)$ is torsion. By the ideal class number formula for
imaginary quadratic field, we have that the special value of
L-series removed infinite factor
$$L^{(\infty)}(1, \eta)=\pi h_K/\sqrt{8n},
\qquad L^{(\infty)} (1, \eta_0)=\pi h_{K_0}/\sqrt{8n_0}.$$

Put all together,  we have that
$$\begin{aligned}
\frac{\wh{h} (P^\chi(f))}{\wh{h} (P^{\chi_0}(f))}&=\frac{L'(1/2, \pi, \chi)}{L'(1/2, \pi, \chi_0)}\cdot \frac{h_{K_0}^{-2}L(1, \eta_0)^2}{h_K^{-2}L(1, \eta)^2}\cdot \prod_{v|\infty 2n}\frac{\beta_v(f)}{\beta_{0, v}(f)}\\
&=\frac{L(1/2, \pi^{(cn/n_0)})}{L(1/2, \pi^{(c)})}\cdot
\frac{n}{n_0}\cdot \prod_{p|(n/n_0)} p^{-1/2}\\
&=\frac{L(1/2, \pi^{(cn/n_0)})/\Omega^{(cn/n_0)}} {L(1/2,
\pi^{(c)})/\Omega^{(c)}}=\frac{L^\alg(E^{(cn/n_0)}, 1)} {L^\alg(
E^{(c)}, 1)},\end{aligned}$$ which completes the proof of
Proposition \ref{Prop B}.
\end{proof}

We have  the following estimation of the  $2$-adic valuation of
$L^\alg(E^{(cn/n_0)}, 1)$. Note that all prime divisor of $n/n_0$
are congruent to $1$ modulo $8$.

\begin{prop}[Zhao \cite{Zhao1}, \cite{Zhao2}] \label{Prop C} Let $s\geq 1$ be an integer and let
$m=p_1\cdots p_s$ be a product of dictinct primes $p_i\equiv 1\mod
8$.  Then
\begin{enumerate}
\item  the $2$-adic additive valuation of  $L^{\alg}(E^{(2m)}, 1)$  is not less than $2s-1$.
\item  the $2$-adic additive valuation of
$L^{\alg}(E^{(m)}, 1)$ is not less than $2s-1$ and the equality
holds if and only if the ideal class group $\CA$ of
$\BQ(\sqrt{-2m})$ satisfies that $\dim_{\BF_2} \CA[4]/\CA[2]=1$ and
$\displaystyle{\prod_{i=1}^s\left(\frac{2}{p_i}\right)_4
\cdot(-1)^{(p_i-1)/8}=-1}$.
\end{enumerate}
\end{prop}
\begin{proof} By Corollary 2 in \cite{Zhao1}, for any $m$ a
product of $s$ distinct odd primes $p_i\equiv 1\mod 4$,  the 2-adic
valuation of $L^\alg(1/2, \pi^{(2m)})$ is not less than $2s-2$  and
is equal to $2s-2$ if and only if there are exactly an odd number of
spanning subtrees in the graph $\wt{G}_{-m}$ whose vertices are $-1,
p_1, \cdots, p_s$ and whose edges are those $(-1, p_i)$ with
$p_i\equiv 5\mod 8$ and those $(p_i, p_j), i\neq j,$ with
$\left(\frac{p_i}{p_j}\right)=-1$. Since all primes $p_i$ in the
proposition are $\equiv 1\mod 8$, the graph $\wt{G}_{-m}$ is not
connected so that the 2-adic valuation can not reach the lower bound
$2s-2$ and therefore (1) follows. The statement (2) is just Theorem
1 in \cite{Zhao2}.

\end{proof}

It is known that  $2$-adic valuations of  $L^{\alg}(1/2, \pi^{(1)})$
and $L^{\alg}(1/2, \pi^{(2)})$ are $-3$ and $-2$, respectively (See
\cite{BSD} p.87). In fact, it is also known that the full BSD
conjecture holds for $E^{(1)}$ and $E^{(2)}$.

Now Theorem \ref{Main Theorem in section 3} follows from
Propositions \ref{Prop A}, \ref{Prop B}, and \ref{Prop C}.

\section{Induction Argument on Quadratic Twists}

We now prove Theorem \ref{Main Theorem 3} by induction on $k$.
Recall that a non-zero integer $m$ is a congruent number if and only
if the Mordell-Weil group $E^{(m)}(\BQ)$ of the elliptic curve
$E^{(m)}: my^2=x^3-x$ has rank greater than zero. Note that
$E^{(m)}\cong E^{(-m)}$ only depends on the square-free part of $m$.
Let $n=p_0p_1\cdots p_k$ be a product of distinct odd primes with
$k\geq 0$, $m=n$ or $2n$ such that $m\equiv 5, 6$, or $7\mod 8$, and
$m^*=(-1)^{(n-1)/2} m$. It is known that $E^{(m)}(\BQ)\cong
E(\BQ(\sqrt{m^*}))^-$,  where $E(\BQ(\sqrt{m^*}))^-$ is the group of
points $P\in E(\BQ(\sqrt{m^*}))$ such that $P^\sigma=-P$ where
$\sigma\in \Gal(\BQ(\sqrt{m^*})/\BQ)$ is the non-trivial element.
Note that the torsion subgroup of $E(\BQ(\sqrt{m^*}))^-$ is $E[2]$.
Thus $m$ is  congruent if we can construct a point $y\in
E(\BQ(\sqrt{m^*}))^-\setminus E[2]$. In this section, we will show
the Heegner divisor $P^\chi(f)\in E(\BQ(\sqrt{m^*}))^-$ defined in
section 3 is of infinite order for $n$ satisfying the condition
(1.1) in Theorem \ref{Main Theorem 2} and the abelian character
$\chi: \Gal(H(i)/K)\ra \{\pm 1\}$ defining the extension
$K(\sqrt{m^*})$ over $K$. In fact, to prove the non-triviality of
$P^\chi(f)$, we will construct a point $y_m\in E(\BQ(\sqrt{m^*}))$
following Monsky \cite{Monsky} which satisfies $4y_m=P^\chi(f)$, and
study its 2-divisibility.

Recall that $K=\BQ(\sqrt{-2n})$ and $\varpi=\sqrt{-2n}\in
K_2^\times$ is a uniformizor at $2$.

\subsection{The Case $n\equiv 1\mod 4$}   We first handle the case with
$p_0\equiv 5\mod 8$.  In this case, $m=n$ and the condition (1.1) in
Theorem \ref{Main Theorem 2} says that the ideal class group $\CA$
of $K=\BQ(\sqrt{-2n})$ has no order $4$ elements, or equivalently
that $2\CA\cong \Gal(H/H_0)$ has odd cardinality.

Let $P\in X_0(32)$ be the image of $i\sqrt{2n}/8$ under the complex
uniformization of $X_0(32)\cong \Gamma_0(32)\bs (\CH \cup
\BP^1(\BQ))$. Let $f: X_0(32)\ra E$ be a degree $2$ modular
parametrization. Define
$$z=f(P)+(1+\sqrt{2}, 2+\sqrt{2})\in E(K^{ab})$$which is actually defined
over $H$ by Theorem \ref{Main Theorem of section 2}, and define
\begin{eqnarray}
y_n=\RTr_{H/K(\sqrt{n})} z\in E(K(\sqrt{n})),\end{eqnarray} which is
our desired point. It turns out that $4y_n=P^\chi(f)$, the point we
defined in Theorem \ref{Main Theorem 3} and studied in section 3.

\begin{thm}\label{Main Theorem A} Let $k\geq 0$ be an integer and $n=p_0p_1\cdots p_k$
a product of distinct primes with $p_0\equiv 5\mod 8$ and $p_1,
\cdots, p_k\equiv 1\mod 8$.  Then the point
$$y_n:=\RTr_{H/K(\sqrt{n})} z$$ is actually defined over
$\BQ(\sqrt{n})$ and  belongs to $2^{k-1} E(\BQ(\sqrt{n}))^-+E[2]$.
Moreover, the point $y_n\notin 2^kE(\BQ(\sqrt{n}))^-+E[2]$ if the
ideal class group of $K=\BQ(\sqrt{-2n})$ does not contain  order $4$
elements.
\end{thm}

We start with the case $k=0$, where the statement $y_n\in
2^{-1}E(\BQ(\sqrt{n}))^-+E[2]$ in the theorem is understood as
$2y_n\in E(\BQ(\sqrt{n}))^-$.  When $k=0$ the above theorem is due
to Monsky \cite{Monsky}:
\begin{prop}Let $n=p_0\equiv 5\mod 8$ be a prime. Then the point $y_{p_0}$
satisfies: $$y_{p_0}\in E(\BQ(\sqrt{p_0}))\setminus
E(\BQ(\sqrt{p_0}))^-\quad\text{and}\quad 2y_{p_0}\in
E(\BQ(\sqrt{p_0}))^-\setminus
\left(2E(\BQ(\sqrt{p_0}))^-+E[2]\right).$$ In particular,
$2y_{p_0}\in E(\BQ(\sqrt{p_0}))^-$ is of infinite order and
therefore $p_0$ is a congruent number.
\end{prop}

\begin{proof}In this case $K(\sqrt{p_0})=H_0$ and the ideal class
group $\CA$ of $K=\BQ(\sqrt{-2p_0})$ satisfies that $2\CA\cong
\Gal(H/H_0)$ has odd cardinality.  Note that the action of the
complex conjugation on $\Gal(H/K)$ is given by inverse and
$\Gal(H/K(\sqrt{p_0}))$ is stable under this action. It follows that
the point $y_{p_0}=\RTr_{H/K(\sqrt{p_0})}z$ is fixed by the action
of complex conjugation and therefore $y_{p_0}\in
E(\BQ(\sqrt{p_0}))$.

For any $t\in \wh{K}^\times$, we have that $z_{\varpi t}+z_t=(0, 0)$
by Theorem \ref{Main Theorem of section 2} (3). Thus
$$y_{p_0}+y_{p_0}^{\sigma_\varpi}=\# \Gal(H/K(\sqrt{p_0})) \cdot (0, 0)=
\# 2\CA\cdot (0, 0)=(0, 0).$$Since
$$\sigma_{\varpi}(\sqrt{p_0})=[(2p_0)_2,
\BQ^{ab}/\BQ](\sqrt{p_0})=[2_2,
\BQ^{ab}/\BQ](\sqrt{p_0})=\left(\frac{2}{p_0}\right)\sqrt{p_0}=-\sqrt{p_0},$$
$\sigma_{\varpi}$ acts on $K(\sqrt{p_0})$ non-trivially,  thus
$2y_{p_0}\in E(\BQ(\sqrt{p_0}))^-$ but $y_{p_0}\notin
E(\BQ(\sqrt{p_0}))^-$.

Suppose that  $2y_{p_0}\in 2E(\BQ(\sqrt{p_0}))^-+E[2]$, say,
$2y_{p_0}=2y'+t$ for some $y'\in E(\BQ(\sqrt{p_0}))^-$ and $t\in
E[2]$. Then we have that $2(y_{p_0}-y')=t$.  Thus $y_{p_0}-y'\in
E[4]\cap E(\BQ(\sqrt{p_0}))=E[2]$ and then $y_{p_0}\in
y'+E[2]\subset E(\BQ(\sqrt{p_0}))^-$, a contradiction. In
particular, $2y_{p_0}\in E(\BQ(\sqrt{p_0}))^-$ is of infinite order
and therefore $p_0$ is a congruent number.
\end{proof}

For each positive divisor $d$ of $n$ divisible by $p_0$, define
$$y_d:=\RTr_{H/K(\sqrt{d})} z\in E(K(\sqrt{d})).$$
Define $y_0=\RTr_{H/H_0} z\in E(H_0)$.
\begin{lem}Assume $k\geq 1$ and let $\varpi\in K_2^\times$ be the uniformizor $\sqrt{-2n}$.
Then \begin{enumerate}
\item the point $y_d\in E(\BQ(\sqrt{d}))^-$ for each $d$ with
$p_0|d|n$.
\item  the point $y_0\in E(H_0^+)$, where $H_0^+=H_0\cap
\BR=\BQ(\sqrt{p_0}, \cdots, \sqrt{p_k})$,  satisfies:
\begin{equation}
y_0+y_0^{\sigma_\varpi}=\# 2\CA \cdot (0, 0).\end{equation}
\end{enumerate}
Moreover these points satisfy the following relation:
\begin{equation}\sum_{p_0|d|n} y_d=\begin{cases}2^k y_0, \qquad
&\text{if $k\geq 2$},\\
2^k y_0+\#2\CA \cdot (0, 0), &\text{if $k=1$.}
\end{cases}\end{equation}
\end{lem}

\begin{proof} Similar to Proposition 4.1, for each $d$ with
$p_0|d|n$, $\Gal(H/K(\sqrt{d}))$ is stable under the action of the
complex conjugation and therefore $y_d\in E(\BQ(\sqrt{d}))$.
Moreover, by Theorem \ref{Main Theorem of section 2},
$$y_d+y_d^{\sigma_\varpi}=\# \Gal(H/K(\sqrt{d}))\cdot (0, 0)=0$$
since the cardinality of $\CA[2]$ is $2^{k+1}$ and thus
$\Gal(H/K(\sqrt{d}))$ has even cardinality $\# 2\CA \cdot 2^k$.
Similarly,  $y_0$ is invariant under the complex conjugation so that
$y_0\in E(H_0^+)$ and satisfies the relation (4.2).

Note that any element $\sigma$ in $\Gal(H_0/K)$ maps $\sqrt{p_i}$ to
$\pm \sqrt{p_i}$ for $0\leq i\leq k$, Note also that $\sigma \in
\Gal(H_0/K(\sqrt{d}))$ if and only if the cardinality of $\{p|d,
\sigma(\sqrt{p})=-\sqrt{p}\}$ is even. For any $1\neq \sigma \in
\Gal(H_0/K)$, let $n_\sigma=\#\{d, p_0|d|n, \sigma \in
\Gal(H_0/K(\sqrt{d}))\}$. Note that
$\sigma_\varpi(\sqrt{p_0})=-\sqrt{p_0}$ and
$\sigma_\varpi\sqrt{p_i}=\sqrt{p_i}$ for all $1\leq i\leq k$.  If
$1\neq \sigma\in \Gal(H/K)$ fixes $\sqrt{p_0}$ and changes exact $s$
elements among $\sqrt{p_1}, \cdots, \sqrt{p_k}$, then
$n_\sigma=2^{k-s} \left[
\binom{s}{0}+\binom{s}{2}+\cdots+\right]=2^{k-1}$ and $n_{\sigma
\sigma_\varpi}=2^{k-s}
\left[\binom{s}{1}+\binom{s}{3}+\cdots\right]=2^{k-1}.$ Thus
$$\begin{aligned}
\sum_{p_0|d|n} y_d-2^k y_0&=\sum_{p_0|d|n} (y_d-y_0)=\sum_{p_0|d|n}
\sum_{1\neq \sigma \in \Gal(H_0/K(\sqrt{d}))} y_0^\sigma\\
&=\sum_{1\neq \sigma \in \Gal(H_0/K)} n_\sigma y_0^\sigma.\\
&=\sum_{1\neq \sigma, \sigma \text{\ fixing $\sqrt{p_0}$}} (n_\sigma y_0^\sigma +n_{\sigma \sigma_\varpi} y_0^{\sigma \sigma_\varpi})\\
&=2^{k-1} (2^k-1) \#2\CA\cdot (0, 0).
\end{aligned}$$
The equality (4.3) now follows.

\end{proof}

\begin{proof}[Proof of Theorem \ref{Main Theorem A}] We do induction on $k$. It holds when $k=0$ by Proposition 4.1.
Now assume  $k\geq 1$.  By Lemma 4.2, we now have
\begin{eqnarray}y_n+\sum_{p_0|d|n, d\neq n} y_d=2^k y_0\mod
E[2].\end{eqnarray} For each $d$ with $p_0|d|n, d\neq n$, let
$y_d^0$ be the point constructed similarly with $\BQ(\sqrt{-2n})$
replaced
 by $\BQ(\sqrt{-2d})$. Then $4y_d, 4y_d^0$ are
$P^\chi(f)$ and $P^{\chi_0}(f)$ in section 3, respectively. Thus
either $y_d^0$ is torsion or the ratio of $y_d$ over $y_d^0$ in the
one dimensional $\BQ$-vector space $E(\BQ(\sqrt{d}))^-\otimes_\BZ
\BQ$ is $[y_d: y_d^0]=[P^\chi(f): P^{\chi_0}(f)]$. By induction
hypothesis and Theorem \ref{Main Theorem in section 3}, we know that
$y_d\in 2^k E(\BQ(\sqrt{d}))^-+E[2]$.

Now write $y_d=2^k y_d'+t_d$ with $y_d'\in E(\BQ(\sqrt{d}))^-$ and
$t_d\in E[2]$. Then we have that
$$y_n=2^k\left(y_0-\sum_{p_0|d|n, d\neq n} y_d'\right)  +t$$
for some $t\in E[2]$. Note any proper sub-extension of $H_0^+/\BQ$
must be ramified at some odd prime. Thus
$$E[2^\infty]\cap E(H_0^+)=E[2].$$

Consider  the (injective) Kummer map
$$E(\BQ(\sqrt{n}))/2^{k+1}E(\BQ(\sqrt{n}))\lra H^1(\BQ(\sqrt{n}), E[2^{k+1}]),$$
and  the exact inflation-restriction sequence
$$1\lra H^1(\Gal(H_0^+/\BQ(\sqrt{n})), E[2])\lra H^1(\BQ(\sqrt{n}), E[2^{k+1}])\lra H^1(H_0^+, E[2^{k+1}]).$$
Since $2y_n=2^{k+1}\left(y_0-\sum_{p_0|d|n, d\neq n} y_d'\right)$
with $y_0-\sum_{p_0|d|n, d\neq n} y_d' \in E(H_0^+)$, we know that
the image of $2y_n$ in the Kummer map belongs to
$H^1(\Gal(H_0^+/\BQ(\sqrt{n})), E[2])$ and then the image is killed
by $2$.  Thus $4y_n\in 2^{k+1}E(\BQ(\sqrt{n}))$, or $4(y_n-2^{k-1}
\wt{y}_n)=0$ for some $\wt{y}_n\in E(\BQ(\sqrt{n}))$. It follows
that $y_n=2^{k-1}\wt{y}_n$ modulo $E[2]$ and then belongs to
$2^{k-1}E(\BQ(\sqrt{n}))+E[2]$. Moreover, with the relation
$$y_n=2^k\left(y_0-\sum_{p_0|d|n, d\neq n} y_d'\right)\mod E[2],$$ we
have that $2^{k-1}\left(\wt{y}_n-2y_0+\sum_{d\neq n} 2y_d'\right)\in
E[2]$,  which implies that $$\wt{y}_n=2y_0-\sum_{d\neq n}
2y_d'+t,\qquad \text{for some $t\in E[2]$}.$$ Note that for any
$0<d|n$ with $p_0|d$, we have that
$\sigma_\varpi(\sqrt{d})=-\sqrt{d}$ and therefore
$y_d^{'\sigma_\varpi}=-y_d'$. Thus
$\wt{y}_n^{\sigma_{\varpi}}+\wt{y}_n=2(y_0^{\sigma_\varpi}+y_0)=0$,
i.e. $\wt{y}_n\in E(\BQ(\sqrt{n}))^-$. This shows that $y_n\in
2^{k-1} E(\BQ(\sqrt{n}))^-+E[2]$.

Now assume that the ideal class group of $K=\BQ(\sqrt{-2n})$ has no
order 4 element. Suppose that $y_n=2^k y_n' +t_n$ for some $y_n' \in
E(\BQ(\sqrt{n}))^-$ and $t_n\in E[2]$. Then by the relation (4.3) in
Lemma 4.2,  we have that
$$2^k\left(y_0-\sum_{p_0|d|n} y_d'\right)\in E[2], \qquad \text{with}\ y_0-\sum_{p_0|d|n} y_d'\in E(H_0^+).$$
Again, any proper sub-extension of $H_0^+/\BQ$ must be ramified at
some odd prime, we have that $y_0-\sum_{p_0|d|n} y_d'=t$ for some
$t\in E[2]$. Thus
$$y_0+y_0^{\sigma_\varpi}=\sum_{p_0|d|n}(y_d'+y_d^{'\sigma_\varpi})+
(t+t^{\sigma_\varpi})=0.$$But $y_0+y_0^{\sigma_\varpi}=(0, 0)\in
E[2]$ by the equality (4.2) in Lemma 4.2.  It is a contradiction.

\end{proof}

\subsection{The Case $n\equiv 3\mod 4$} In this subsection, we assume that
$n=p_0p_1\cdots p_k$ is a product of distinct primes with $p_0\equiv
3\mod 4$ and $p_1, \cdots, p_k\equiv 1\mod 8$.

The genus field  of $K=\BQ(\sqrt{-2n})$ is $H_0=K(\sqrt{-p_0},
\sqrt{p_1}, \cdots, \sqrt{p_k})$. Thus $\sqrt{2}\in H_0$ but
$i\notin H$.  We identify the ideal class group $\CA$ of $K$ with
the subquotient of $\wh{K}^\times$ corresponding to
$\Gal(H(i)/K(i))$. Recall that we denote by $\varpi$ the uniformizer
$\sqrt{-2n}$ in $K_2^\times$ and let $\varpi'=\varpi (1+\varpi)\in
K_2^\times$ which is also a uniformizer and $\sigma_{\varpi'}$ fixes
$i$. Let $\varpi_{p_i}'=\sqrt{-2n}\in K_{p_i}^\times$ for $1\leq
i\leq k$ and let $\varpi_{p_0}'=(\sqrt{-2n})_{p_0} (1+\varpi)\in
K_{p_0}^\times K_2^\times$. Note that the condition (1.1):
$\dim_{\BF_2} \CA[4]/\CA[2]=1$,  in Theorem \ref{Main Theorem 2}
says that
\begin{itemize}
\item if $p_0\equiv 3\mod 8$, then $\CA$ has no order 4 elements, or
equivalently, $2\CA$ has odd cardinality.
\item if $p_0\equiv 7\mod 8$, then the class of $\varpi'\in K_2^\times$
in $\CA$ is the only non-trivial element in $\CA[2]\cap 2\CA$. In
fact, by Gauss' genus theory, one can check that
$\sigma_{\varpi'}\big|_H$ fixes $\sqrt{-p_0}$ and all $\sqrt{p_i}$
for $1\leq i\leq k$.
\end{itemize}

Let $d\equiv 6, 7\mod 8$ be a  positive divisor of $2n$, then
$\sqrt{-d}\in H_0$ is fixed under $\sigma_{\varpi'}$.  Let
$\chi=\chi_d$ be the character of $\CA=\Gal(H(i)/K(i))$ factoring
through $\Gal(K(i, \sqrt{-d})/K(i))$ which is non-trivial when
$d\neq 2n$.  Since $\Ker \chi$ contains the class $[\varpi']\in \CA$
of $\varpi'$, the character $\chi$ factors through $\CA/\langle
[\varpi']\rangle$. For any complete representatives $\phi\subset
\CA$ of $\CA/[\varpi']$, we define \begin{eqnarray}y_{d,
\phi}=\sum_{t\in \phi} \chi(t) z_t.\end{eqnarray} The point $y_{d,
\phi}$ is indepent of $\phi$ up to $E[2]$ by Theorem \ref{Main
Theorem 2 in section 2} (3). More precisely,  let $\phi'$ be a
second set of representatives and $\phi^c$ the complement of $\phi$
in $\CA$. If $\phi'\cap \phi^c$ has even cardinality then $y_{d,
\phi'}=y_{d, \phi}$; and if the cardinality is odd, then $y_{d,
\phi'}-y_{d, \phi}=(-1, 0)$ or $(1, 0)$ according as $n\equiv 3$ or
$7\mod 8$. We will often ignore the dependence of $y_{d, \phi}$ on
$\phi$ and abbreviate it to $y_d$. Note that $4y_d=P^\chi(f)$ with
$m_0=d$ in section 3 in then case $n\equiv 3\mod 4$.

Recall that $E(\BQ(\sqrt{-d}))^-$ denotes the subgroup of points
$P\in E(\BQ(\sqrt{-d}))$ satisfying $\sigma P=-P$ where $\sigma$ is
the non-trivial  element in  $\Gal(\BQ(\sqrt{-d})/\BQ)$. Then
$E(\BQ(\sqrt{-d}))^-$ is isomorphic to the Mordell-Weil group of the
elliptic curve $E^{(d)}: dy^2=x^3-x$ over $\BQ$ and its torsion
subgroup is $E[2]$.

Let $m=n$ or $2n$ such that $m\equiv 6$ or $7\mod 8$. We will show
that $y_m$ is of infinite order under the condition (1.1) in Theorem
\ref{Main Theorem 2} by induction on the number of prime divisors of
$n$. More precisely, the remain of this section will be devoted to
prove the following two theorems.

\begin{thm}\label{Main Theorem B} Let $k\geq 0$ an integer
and $n=p_0p_1\cdots p_k$ a product of distinct primes with
$p_0\equiv 7\mod 8$ and $p_1, \cdots, p_k\equiv 1\mod 8$. Then for
$m=n$ or $2n$, $y_m\in 2^{k-1} E(\BQ(\sqrt{-m}))^-+E[2]$, and if
$\dim_{\BF_2}\CA[4]/\CA[2]=1$ then the point $y_m\notin
2^kE(\BQ(\sqrt{-m}))^-+E[2]$.
\end{thm}

\begin{thm}\label{Main Theorem C} Let $k\geq 0$
be an integer and $n=p_0p_1\cdots p_k$ a product of distinct primes
with $p_0\equiv 3\mod 8$ and $p_1, \cdots, p_k\equiv 1\mod 8$. Then
$y_{2n}\in 2^{k-1} E(\BQ(\sqrt{-2n}))^-+E[2]$,   and if the field
$K=\BQ(\sqrt{-2n})$ has no order 4 ideal class, then the point
$y_{2n}\notin 2^kE(\BQ(\sqrt{-2n}))^-+E[2]$.
\end{thm}
The following proposition is the initial case for the induction
process and proved in \cite{Monsky}, we repeat its proof here for
completeness.

\begin{prop} \label{Initial 2} Let $n=p_0$ be a prime congruent to $3$ modulo $4$.
Let $m=p_0$ or $2p_0$ such that $m\equiv 6, 7\mod 8$.  Then we have
$$
2y_m\in E(\BQ(\sqrt{-m}))^-\setminus
\left(2E(\BQ(\sqrt{-m}))^-+E[2]\right).$$ In particular, $2y_m\in
E(\BQ(\sqrt{-m}))^-$ is of infinite order and therefore $m$ is a
congruent number.
\end{prop}

\begin{proof} Note that the Galois group $\Gal(H(i)/\BQ)$ is generated by
$\Gal(H(i)/K(i))$, the complex conjugation, and the operator
$\sigma_{1+\varpi}$.

The element $\sigma_{1+\varpi}$ induces the non-trivial involution
on $H(i)$ over $H$.  By Theorem \ref{Main Theorem 2 in section 2}
(1),
\begin{eqnarray} y_{m,
\phi}^{\sigma_{1+\varpi}}-y_{m, \phi}=\sum_{\phi} (0,
0)\end{eqnarray} Thus $\sigma_{1+\varpi}$ fixes $2y_m$, i.e.  $2y_m$
is rational over $H$. Moreover, by Theorem \ref{Main Theorem 2 in
section 2} (2),
\begin{eqnarray}\ov{y_{m, \phi}}=\sum_{[t]\in \phi} \chi(t)
(-z_{t^{-1}}+(1, 0))=-y_{m, \phi^{-1}}+\sum_\phi (1,
0),\end{eqnarray} thus the complex conjugation of $2y_m$ is equal to
$-2y_m$. Any $\sigma_s \in \Gal(H(i)/K(i))$ maps $y_{m, \phi}$ to
\begin{eqnarray}\sum_\phi \chi(t) z_{st}=\chi(s) y_{m,
[s]\phi}.\end{eqnarray} So it takes $2y_m$ to $2y_m$ or $-2y_m$
according to $\sigma_s$ acts trivially or non-trivially on
$\sqrt{-m}$. Thus $2y_m$ is rational over $K(i, \sqrt{-m})$ and
therefore rational over $K(i, \sqrt{-m})\cap H=K(\sqrt{-m})$. Now we
claim that $2y_m$ is rational over $\BQ(\sqrt{-m})$.  This is clear
if $m=2n$. When $m\neq 2n$ choose an element $\sigma\in
\Gal(H(i)/K(i))$ mapping $\sqrt{-m}$ to $-\sqrt{-m}$, then both
$\sigma$ and the complex conjugation take $2y_m$ to $-2y_m$, and
therefore their composition fixes $2y_m$ and has fixed field
$\BQ(\sqrt{-m})$ in $K(\sqrt{-m})$. This shows the claim and
therefore it follows that $2y_m\in E(\BQ(\sqrt{-m}))^-$.

Let us first consider the case with  $p_0\equiv 3\mod 8$. Then
$m=2p_0\equiv 6\mod 8$ and $\phi$ has odd cardinality. We need to
show that $2y_m\notin 2E(\BQ(\sqrt{-m}))^-+E[2]$. Suppose this is
not the case, i.e. $2y_m=2y+t$ for some $y\in E(\BQ(\sqrt{-m}))^-$
and $t\in E[2]$. Then $P:=y_m-y$ is a 4-torsion point. Note that
$\sqrt{-m}\in H$. Then, by the equation (4.4), we have that
$$\sigma_{1+\varpi} (P)-P=\sigma_{1+\varpi}(y_m)-y_m
=\sum_\phi (0, 0)=(0, 0).$$ On the other hand,  as we pointed out at
the beginning of section 2, $E[4]/E[2]$ is represented by $0, (i,
1-i), (1+\sqrt{2}, 2+\sqrt{2}), (-1-\sqrt{2}, i(2+\sqrt{2}))$. Note
that $\sigma_{1+\varpi}$ moves $i$ but fixes $\sqrt{2}$. Thus
$\sigma_{1+\varpi}$  acts on any point $Q\in E[4]$ via the complex
conjugation. It follows immediately that  $\sigma_{1+\varpi}
(Q)-Q=0$ if $Q\equiv 0, (1+\sqrt{2}, 2+\sqrt{2})\mod E[2]$ and
$\sigma_{1+\varpi}(Q)-Q=(-1, 0)$ otherwise. It is a contradiction.

Assume now that  $p_0\equiv 7\mod 8$. Then $m=p_0$ or $2p_0$
corresponding $\chi$ non-trivial or trivial.  We need to show that
$2y_m\notin 2E(\BQ(\sqrt{-m}))^-+E[2]$. Suppose this is not the
case, i.e. $2y_m=2y\mod E[2]$ for some $y\in E(\BQ(\sqrt{-m}))^-$.
Since $\phi$ has even cardinality,  $y_m$ is rational over $H$. Then
$P:=y_m-y\in E[4]\cap E(H)=E[4]\cap E((\BQ(\sqrt{2}))$ and therefore
$P=0$ or $(1+\sqrt{2}, 2+\sqrt{2})$ modulo $E[2]$.

Note that $\CA[2^\infty]$ is cyclic by Gauss' genus theory. Take an
element $[t]\in\CA-2\CA$ (for example, a generator of
$\CA[2^\infty]$), then $\CA/([t])$ has odd cardinality. Let $\phi_0$
be a set of representatives for the $\CA/([t])$ and then we may take
$\phi=\bigcup_{i=0}^{m-1} [t]^i \phi_0$ if the order of $[t]$ is
$2m$ (thus $[t]^m=[\varpi']$).  Use this $\phi$ to define $y_m$. Now
we have
$$\begin{aligned}
y_{m, \phi}^{\sigma_t}&=\left(\sum_{[t']\in\phi}
\chi(t')z_{t'}\right)^{\sigma_t}
=\chi (t)  \sum_{[t']\in [t]\phi}\chi(t')z_{t'}\\
&=\chi(t) \left( \sum_\phi \chi(t')
z_{t'}+\sum_{\phi_0}\chi(t')(z_{\varpi't'}-z_{t'})\right)=\chi(t)y_{m,
\phi}+(1, 0)
\end{aligned}$$
Note that $\sigma_t$  fixes $\sqrt{-2p_0}$ but moves $\sqrt{2}$ and
$\sqrt{-p_0}$. Thus $\sigma_t y=\chi(t)y$ and then
$$P^{\sigma_t}-\chi(t)P=y_m^{\sigma_t}-\chi(t) y_m=(1, 0).$$But we have shown that
$P=0$ or $(1+\sqrt{2}, 2+\sqrt{2})$ modulo $E[2]$. If follows that
$P^{\sigma_t}-\chi(t)P=0$ if $P=0\mod E[2]$ and that
$P^{\sigma_t}-\chi(t)P=(-1, 0)$ or $0, 0)$ if $P=(1+\sqrt{2},
2+\sqrt{2})\mod E[2]$ according to $\chi(t)=1$ or $-1$. It is a
contradiction.
\end{proof}
We now refine the beginning argument of previous Proposition to
obtain the defining field of points $y_d$'s when $k\geq 1$.
\begin{lem} Assume that $k\geq 1$ and that $n=p_0p_1\cdots p_k$ is a product
of distinct primes $p_0\equiv 3\mod 4$ and $p_i\equiv 1\mod 8$ for
$1\leq i\leq k$.  Let $d|2n$ be a positive integer $\equiv 6, 7\mod
8$. Then the point $y_d\in E(\BQ(\sqrt{-d}))^-$.
\end{lem}
\begin{proof}Note again that the Galois group $\Gal(H(i)/\BQ)$ is generated by
$\Gal(H(i)/K(i))$, the complex conjugation, and the operator
$\sigma_{1+\varpi}$. Take a prime $p|n$, then the class
$[\varpi_p']\in \CA[2]$ is neither trivial nor equal to  $[\varpi']$
since $k\geq 1$. Let $\phi_0\subset \CA$ be a complete set of
representatives of $\CA/([\varpi'], [\varpi_p'])$ and let
$\phi=\phi_0\cup [\varpi_p'] \phi_0$, which has even cardinality.
\begin{itemize}
\item[(i)] Since $\phi$ has even cardinality, by Theorem \ref{Main Theorem 2 in section
2} (1), we have
 $$y_{d, \phi}^{\sigma_{1+\varpi}}-y_{d, \phi}=\sum_{\phi} (0,
0)=0.$$ Thus $\sigma_{1+\varpi}$ fixes $y_d$, i.e.  $y_d$ is
rational over $H$.
\item[(ii)] The set $\phi^{c}\cap \phi^{-1}$ is stable under multiplication by
$[\varpi_p']$ and then has even cardinality. It follows that
$y_{d,\phi}=y_{d, \phi^{-1}}$. Then by Theorem \ref{Main Theorem 2
in section 2} (2), we have
$$\ov{y_{d, \phi}}=\sum_{[t]\in \phi}
\chi(t) (-z_{t^{-1}}+(1, 0))=-y_{d, \phi^{-1}}+\sum_\phi (1,
0)=-y_{d, \phi}.$$Thus the complex conjugation maps $y_d$ to its
negative.
\item[(iii)]  For any $[t]\in \CA$,
$[t]\phi=[t]\phi_0\cup [t\varpi_p']\phi_0$, the set $\phi^{c}\cap
[t]\phi$ is stable under multiplication by $[\varpi_p']$ too and
then has even  cardinality. Thus for any $\sigma_t \in
\Gal(H(i)/K(i))$,
$$y_{d, \phi}^{\sigma_t}=\sum_\phi
\chi(t') z_{tt'}=\chi(t) y_{d, [t]\phi}=\chi(t)y_{d, \phi},$$ So
$y_d^{\sigma_t}$ is equal to $y_d$ or $-y_d$ according to $\sigma_t$
acts trivially or non-trivially on $\sqrt{-d}$, thus $y_d$ is
rational over $K(i, \sqrt{-d})$.
\end{itemize}
Now by the same argument as the previous Proposition, we have that
the point $y_d$ belong to $E(\BQ(\sqrt{-d}))^-$.
\end{proof}

We are now going to show separately Theorem \ref{Main Theorem B} and
Theorem \ref{Main Theorem C}, which correspond to the case with
$p_0\equiv 7\mod 8$ and the case with $p_0\equiv 3\mod 8$.

Let $k\geq 1$ be an integer and $n=p_0p_1\cdots p_k$ with $p_0\equiv
7\mod 8$ and $p_i\equiv 1\mod 8$ for $1\leq i\leq k$. Note that
$[\varpi']\in 2\CA\cap \CA[2]$.  Let $\phi_0$ be a set of
representatives of $2\CA/([\varpi'])$. Let $\psi$ be a set of
representatives of $\CA/2\CA$. Then $\phi=\bigcup_{[s]\in \psi}
[s]\phi_0$ is a set of representatives of $\CA/([\varpi'])$. We use
this $\phi$ to define all $y_d$'s. Let $\beta\in \Gal(H(i)/K(i))$ be
an element which moves $\sqrt{2}, \sqrt{p_1}, \cdots, \sqrt{p_k}$.
Then  $\beta$ fixes or moves $\sqrt{p_0}$ according to that $k$ is
odd or even.

\begin{lem} Assume $k\geq 1$. Let $n=p_0p_1\cdots p_k$ with
$p_0\equiv 7\mod 8$ and $p_i\equiv 1\mod 8$ for $1\leq i\leq k$. Let
$m=n$ or $2n$ such that $m\equiv 6$ or $7$ modulo $8$. Then we have
\begin{enumerate} \item for each positive divisor $d$ of $2n$ congruent to $6$
or $7\mod 8$,  the point $y_d\in E(\BQ(\sqrt{-d}))^-$.
\item Let $y_0:=\sum_{[t]\in \phi_0} z_t$. Then $y_0+ (-1)^m y_0^\beta$ is
rational over the genus subfield $H_0=K(\sqrt{-p_0},
\sqrt{p_1},\cdots, \sqrt{p_k})$ of $K$.
\end{enumerate}
These points satisfy the following relation:
\begin{equation}\sum_{\substack{{p_0|d|2n}\\ \nu_0(d)\equiv\nu_0(m)\mod 2}}
y_d=2^k\left(y_0+(-1)^m y_0^\beta \right)\end{equation} where for
any integer $d$,  $\nu_0(d)$ denotes the number of prime divisors of
$d$.
\end{lem}

\begin{proof}Note that the Galois group $\Gal(H(i)/H_0)$ is generated by $\Gal(H(i)/H_0(i))$
and the operator $\sigma_{1+\varpi}$. The statement (1) is showed in
Lemma 4.7. By Theorem \ref{Main Theorem 2 in section 2}, for any
$[s]\in 2\CA$, $y_0-y_0^{\sigma_s}=0$ or $(1, 0)$ and thus $(y_0\pm
y_0^\beta)$ is fixed by $\sigma_s$ and then rational over $H_0(i)$.
Since $\sigma_{1+\varpi}$ induces the involution of $H(i)$ over $H$,
Theorem \ref{Main Theorem 2 in section 2} also implies that $y_0\pm
y_0^\beta$ is rational over $H_0$. Moreover, we  have that
$$\sum_{\substack{{p_0|d|2n}\\ \nu_0(d)\equiv\nu_0(m)\mod 2}}
y_d=\sum_{[t]\in
\phi_0}\sum_{[s]\in\psi}\left(\sum_{\substack{{p_0|d|2n}\\
\nu_0(d)\equiv\nu_0(m)\mod 2}} \chi_d(s)\right) (z_t)^{\sigma_s}$$
It is clear that the summation in the last bracket is equal to $2^k$
for $\sigma_s=1$, $(-1)^m\cdot 2^k$ for $\sigma_s=\beta$, and $0$
otherwise. Thus the equality (4.9) follows .
\end{proof}

\begin{proof}[Proof of Theorem \ref{Main Theorem B}] We prove the theorem by induction on $k$.
 The initial case $k=0$ is given by Proposition \ref{Initial 2}. Now assume that $k\geq 1$.  For $m=n$ or $2n$,
Similar to the case with $p_0\equiv 5\mod 8$, we have the following
\begin{enumerate}
\item the point $y_m\in 2^{k-1} E(\BQ(\sqrt{-m}))^-+E[2]$, using the equality (4.9);
\item for each positive $d$ with $p_0|d|2n$ and $d\neq n, 2n$, $y_d\in 2^k E(\BQ(\sqrt{-d}))^-+E[2]$, i.e. of form $2^k y_d'+t_d$ for some $y_d'\in E(\BQ(\sqrt{-d}))^-$ and $t_d\in E[2]$.
\end{enumerate}
Now we show that $y_m\notin 2^k E(\BQ(\sqrt{-m}))^-+E[2]$ under the
condition $\dim_{\BF_2}\CA[4]/\CA[2]=1$. Suppose it is not the case,
i.e. $y_m=2^k y_m'+t_m$ for some $y_m'\in 2^kE(\BQ(\sqrt{-m}))^-$
and $t_m\in E[2]$. Thus the previous lemma implies that
$$P:=\left(y_0+(-1)^{m} y_0^\beta -\sum_{\substack{{p_0|d|2n}\\ \nu_0(d)\equiv\nu_0(m)\mod 2}} y_d'\right)\in E[2^{k+1}]\cap E(H_0)=
E[4]\cap E(\BQ(\sqrt{2})).$$ It follows that $P^\beta-P=0$ or
$(0,0)$ and $P^\beta+P=0$ or $(-1, 0)$.

The assumption $\dim_{\BF_2} \CA[4]/\CA[2]=1$ implies that
$[\varpi']$ is the unique non-trivial element in $2\CA\cap \CA[2]$.
Note that $(y_d')^\beta=(-1)^m y_d'$ and therefore
$$P-(-1)^m P^\beta=y_0-y_0^{\beta^2}.$$
Write $\beta=\sigma_ {t_0}$, then $[t_0]\in \CA\setminus 2\CA$ has
order $4s$ for some integer $s$. Then $[t_0]^{2s}=[\varpi']$ and the
group $2\CA/([t_0]^2)$ is of odd order. Let $\phi_1$ be a set of
representatives for the group $2\CA/([t_0]^2)$, then we may take
$\displaystyle{\phi_0=\bigcup_{i=0}^{s-1}[t_0]^{2i}\phi_1}$ to be
our set of representative for $2\CA/([\varpi'])$ and use it to
define $y_0$. Then
$$
y_0-y_0^{\beta^2}=\sum_{[t]\in \phi_0}z_t-\sum_{[t]\in
\phi_0}z_{t_0^2t}=\sum_{[t]\in\phi_1} (z_t-z_{\varpi' t})=\#
2\CA/([t_0]^2)\cdot (1, 0)=(1, 0).$$It is a contradiction.
\end{proof}
\bigskip

Finally, we  consider the case that $p_0\equiv 3\mod 8$ and $k\geq
1$. Then $\CA=2\CA \times \CA[2]$ and $[\varpi']\in \CA[2]$. Let
$\psi$ be the set of representatives for $\CA[2]/([\varpi'])$
consisting of those $[\varpi'_d]$ fixing $\sqrt{2}$, then
$\phi:=\bigcup_{[s]\in \psi} [s] (2\CA)$ is  a set of
representatives for $\CA/([\varpi'])$. The set $\phi$ is stable
under $[t]\mapsto [t]^{-1}$ and  we use $\phi$ to define all
$y_d$'s.

\begin{lem} Assume $k\geq 1$. Let $n=p_0p_1\cdots p_k$ be a product of distinct
primes with  $p_0\equiv 3\mod 8$ and $p_i\equiv 1\mod 8$ for $1\leq
i\leq k$. Then we have
\begin{enumerate}
\item for each positive divisor $d$ of $2n$ divisible by $2p_0$,
the point $y_d\in E(\BQ(\sqrt{-d}))^-$;
\item the point $y_0:=\sum_{[t]\in 2\CA} z_t\in E(H_0(i))$ satisfies the following
relation:
\begin{equation}
y_0^{\sigma_{1+\varpi}}-y_0=\# 2\CA \cdot (0, 0)
\end{equation}\end{enumerate}Moreover, these points satisfy the following relation:
\begin{equation}\sum_{2p_0|d|2n} y_d=2^ky_0.\end{equation}
\end{lem}
\begin{proof}
The statement (1) is showed in Lemma 4.7. The statement (2) follows
from Theorem \ref{Main Theorem 2 in section 2} (1):
$$y_0^{\sigma_{1+\varpi}}-y_0=\sum_{[t]\in 2\CA}
(z_t^{\sigma_{1+\varpi}}-z_t)=\# 2\CA \cdot (0, 0).$$ Moreover,
$$\sum_{2p_0|d|2n} y_d=\sum_{[t]\in 2\CA} \sum_{[s]\in \psi} \left( \sum_{2p_0|d|2n}
\chi_d(s) \right)z_t^{\sigma_s}.$$ It is clear that the summation in
the last bracket is equal to $2^k$ for $\sigma_s=1$ and $0$
otherwise. Thus the equality (4.11) follows.
\end{proof}

\begin{proof}[Proof of Theorem \ref{Main Theorem C}] We prove the theorem by induction on $k$.
The initial case $k=0$ is given by Proposition \ref{Initial 2}. Now
assume that $k\geq 1$. Similar to previous cases, we have that
\begin{enumerate}
\item the point $y_{2n}\in 2^{k-1} E(\BQ(\sqrt{-2n}))^-+E[2]$, using the equality (4.11);
\item for each positive $d$ with $2p_0|d|2n$ and $d\neq 2n$, $y_d\in 2^k E(\BQ(\sqrt{-d}))^-+E[2]$, i.e. of form $2^k y_d'+t_d$ for some $y_d'\in E(\BQ(\sqrt{-d}))^-$ and $t_d\in E[2]$.
\end{enumerate}
Now we show that $y_{2n}\notin 2^k E(\BQ(\sqrt{-2n}))^-+E[2]$ if the
field $K=\BQ(\sqrt{-2n})$ has no order 4 ideal class. Suppose it is
not the case, i.e. $y_{2n}=2^k y_{2n}'+t_{2n}$ for some $y_{2n}'\in
2^kE(\BQ(\sqrt{-2n}))^-$ and $t_{2n}\in E[2]$. Then as before, we
have that
$$P:=y_{2n}'-y_0+\sum_{2p_0|d|2n, d\neq 2n} y_d'\in E[2^{k+1}]\cap E(H_0(i))=E[2^{k+1}]\cap E(\BQ(i, \sqrt{2}))=E[4].$$
Thus we have the formula
$$y_0=\sum_{2p_0|d|2n} y_d'-P$$
with $P\in E[4]$, and then
$$y_0^{\sigma_{1+\varpi}}-y_0=\sum_{2p_0|d|2n}
((y_d')^{\sigma_{1+\varpi}}-y_d')-(P^{\sigma_{1+\varpi}}-P)=-\ov{P}+P=0\quad
\text{or}\quad (-1, 0).$$ But if $\CA$ has no order 4 element or
equivalently $2\CA$ is odd, we have $y_0^{\sigma_{1+\varpi}}-y_0=(0,
0)$ by the equality (4.10) in Lemma 4.9. It is a contradiction.

\end{proof}

\section{Proof of Main Results}
In this section, we give proofs of Theorem \ref{Easy Theorem},
Theorem \ref{Main Theorem 2}, and Theorem \ref{Main Theorem 3}.

\begin{proof}[Proof of Theorem \ref{Main Theorem 3}] Let
$n=p_0p_1\cdots p_k$ be a product of distinct odd primes with
$p_i\equiv 1\mod 8$, $1\leq i\leq k$, and $p_0\nequiv 1\mod 8$. Let
$m=n$ or $2n$ such that $m\equiv 5, 6$ or $7$ modulo $8$, and
$m^*=(-1)^{(n-1)/2}m$. Let $\chi$ be the abelian character over $K$
defining the unramified extension $K(\sqrt{m^*})$.  Let $P\in
X_0(32)$ be the point $[i\sqrt{2n}/8]$ if $n\equiv 5\mod 8$ and
$[(i\sqrt{2n}+2)/8]$ if $n\equiv 6, 7\mod 8$. Then by Theorem
\ref{Main Theorem of section 2} and Theorem \ref{Main Theorem 2 in
section 2}, we have $f(P)\in E(H(i))$. We showed in Lemma 3.2 that
the point
$$P^\chi(f)=\sum_{\sigma \in \Gal(H(i)/K)}f(P)^\sigma \chi(\sigma)$$
belongs to $E(\BQ(\sqrt{m^*}))^-$ by taking $m_0=m$ there.  It is
easy to see from the definitions (4.1) and (4.5) of $y_m$ that
$$P^\chi(f)=4y_m.$$
By Theorem \ref{Main Theorem A}, Theorem \ref{Main Theorem B}, and
Theorem \ref{Main Theorem C}, we know that for integers $m$ in
Theorem \ref{Main Theorem 3},
$$y_m\in 2^{k-1}E(\BQ(\sqrt{m^*}))^-+E[2]\setminus
2^{k}E(\BQ(\sqrt{m^*}))^-+E[2].$$It then follows that, by noting
that $E[2^\infty]\cap E(\BQ(\sqrt{m^*}))^-=E[2]$,
$$P^\chi(f)\in 2^{k+1} E(\BQ(\sqrt{m^*}))^-\setminus
\left(2^{k+2}E(\BQ(\sqrt{m^*}))^-+E[2]\right).$$In particular,
$P^\chi(f)\in E(\BQ(\sqrt{m^*}))^-$ is of infinite order and $m$ is
a congruent number. This completes the proof of our main result
Theorem \ref{Main Theorem 3}.
\end{proof}

By the following lemma, the condition (1.1) in Theorem \ref{Main
Theorem 2} is actually easy to check.  This allow us not only to
show the existence result Theorem \ref{Easy Theorem} but also to
construct many congruent numbers.
\begin{lem} Let $n=p_0p_1\cdots p_k$ be a product of distinct odd
primes with $p_i\equiv 1\mod 8$ for $1\leq i\leq k$ and $p_0\nequiv
1\mod 8$. Let $\CG$ be the graph whose vertices set $V$ consists of
$p_0, \cdots, p_k$ and whose edges are those $p_ip_j$, $i\neq j$,
with the quadratic residue symbol $\left(\frac{p_i}{p_j}\right)=-1$.
Then the condition (1.1) in Theorem \ref{Main Theorem 2} is
equivalent to any one of the following conditions:
\begin{enumerate}
\item  there does not exist a proper even partition of vertices
$V=V_0\cup V_1$ in the sense that any $v\in V_i$ has even number
edges to $V_{1-i}$, $i=0, 1$; \item the graph $\CG$ has exactly odd
number of spanning subtrees.\end{enumerate}
\end{lem}
\begin{proof} Note that the multiplication by $2$ induces an isomorphism
$\CA[4]/\CA[2]\stackrel{\times 2}{\simeq} \CA[2]\cap 2\CA$.  The
condition (1.1) is the same as
$$2\CA\cap \CA[2]=\begin{cases} 0,
\quad &\text{if $n\equiv \pm 3\mod 8$};\\
 \{0, [\varpi]\}, &\text{otherwise.}\end{cases}$$ Note that the
group $\CA[2]$ consists of $[\varpi_d]$ for all positive divisors
$d|n$, that  $[\varpi]=[\varpi_n]$ in $\CA$, and  that
$[\varpi_d]\in 2\CA$ if and only if
$$\left(\frac{d}{p}\right)=1, \forall\ p|(n/d) \quad \text{and}\quad
\left(\frac{2n/d}{p}\right)=1, \forall\ p|d.$$ The equivalence
between (1) and the condition (1.1) in Theorem \ref{Main Theorem 2}
is then clear. See either \cite{Feng} Lemma 2.2. or \cite{Li-Tian}
Lemma 2 for the equivalence between (1) and (2).

\end{proof}
By Dirichlet theorem, we can construct infinitely many numbers $n$
for each given isomorphism class of graph with exact odd number
spanning subtrees. Thus Theorem \ref{Easy Theorem} follows from
Theorem \ref{Main Theorem 3}. We can even obtain the following
stronger version.
\begin{thm}Let $p_0\nequiv 1\mod 8$ be an odd prime. Then there
exists an infinite set $\Sigma$ of primes congruent to $1$ modulo
$8$ such that the product of $p_0$ (resp. $2p_0$) and primes in any
finite subset of $\Sigma$ is a congruent number if $p_0\equiv 5,
7\mod 8$ (resp. $p_0\equiv 3\mod 4$).
\end{thm}
\begin{proof}
Suppose we are given an odd prime $p_0\nequiv 1\mod 8$. By Dirichlet
theorem, we can choose inductively primes $p_1, p_2, \cdots$
satisfying the following conditions: \begin{itemize} \item[(i)] all
$p_i\equiv 1\mod 8$, \item[(ii)] the quadratic residue symbol
$\left(\frac{p_i}{p_0}\right)=-1$, and
\item[(iii)] for all $1\leq j\leq i-1$,  $\left(\frac{p_i}{p_j}\right)=1$.\end{itemize} Let
$\Sigma$ be the infinite set $\{p_1, p_2, \cdots\}$. Then the graph
with vertexes set $\{p_0\}\cup \Sigma$ and edges those $p_ip_j$
satisfying $\left(\frac{p_i}{p_j}\right)=-1$,  is an infinite
star-shape graph. It is now easy to see that the product $n$ of
$p_0$ with primes in any finite subset of $\Sigma$ satisfies the
condition (1.1) in Theorem \ref{Main Theorem 2}. Thus the set
$\Sigma$ is the desired.
\end{proof}

\medskip

Finally, we explain how Theorem \ref{Main Theorem 2} follows from
Theorem \ref{Main Theorem 3}. Recall $E'=(X_0(32), [\infty])$ is an
elliptic curve defined over $\BQ$ with Weierstrass equation
$y^2=x^3+4x$.
\begin{lem} Let $n=p_0p_1\cdots p_k$ be a product of distinct odd
primes with $p_i\equiv 1\mod 8$ for $1\leq i\leq k$. Let $m=n$ or
$2n$ such that $m\equiv 5, 6$, or $7\mod 8$. Let $\varphi$ be a
2-isogeny from $E^{(m)}$ to $E^{'(m)}: my^2=x^3+4x$ with kernel $0,
(0, 0)$, and $\psi$ its dual isogeny. Let $S^{(\varphi)}$ (resp.
$S^{(\psi)}$) denote the $\varphi$-(resp. $\psi$-) Selmer group.
Then the condition (1.1) in Theorem \ref{Main Theorem 2} implies
that
$$\dim_{\BF_2} S^{(\varphi)}=1\quad\text{and}\quad \dim_{\BF_2}
S^{(\psi)}=2,$$and moreover, that  the $2$-Selmer group
$S^{(2)}(E^{(m)}/\BQ)$ satisfies
\begin{eqnarray}\dim_{\BF_2}
S^{(2)}(E^{(m)}/\BQ)\Big/E^{(m)}[2]=1.\end{eqnarray}
\end{lem}
\begin{proof}
These Selmer groups can be computed using \cite{Silverman}  Ch.X
Proposition 4.9. For example, when $m=n\equiv 5\mod 8$, the
$\varphi$-Selmer group $S^{(\varphi)}\subset H^1(\BQ,
E^{(m)}[\varphi])\cong \BQ^\times/\BQ^{\times 2}$ have
representatives divisors $d$ (including negative ones) of $2n$
satisfying that the curve
$$C_d:\quad dw^2=d^2+ 4n^2z^4$$ over $\BQ$ is solvable locally
everywhere. the $\psi$-Selmer group $S^{(\psi)}\subset H^1(\BQ,
E^{'(m)}[\psi])\cong \BQ^\times/\BQ^{\times 2}$ have representatives
divisors $d$ (including negative ones) of $2n$ satisfying that the
curve
$$C_d':\quad dw^2=d^2-n^2z^4$$ over $\BQ$ is solvable locally
everywhere. Using the equivalence between the condition (1.1) in
Theorem \ref{Main Theorem 2} and (1) in Lemma 5.1, one can check
easily by Hensel Lemma that $S^{(\varphi)}$ consists of only two
elements $1, n$ and is then of $\BF_2$-dimension 1, and $S^{(\psi)}$
consists of $\pm 1, \pm n$ and then is of dimension $2$.

Note that there is an exact sequence:
$$0\lra S^{(\varphi)} \lra S^{(2)} (E^{(m)}/\BQ) \lra S^{(\psi)}.$$
Since the subgroup $E^{(m)}[2]\subset S^{(2)}(E^{(m)}/\BQ)$ of
2-torsion points  provides 2-dimensional image  in $S^{(\psi)}$, the
last morphism is also surjective and therefore $\dim_{\BF_2}
S^{(2)}(E^{(m)}/\BQ)=3$.  The formula (5.1) follows.
\end{proof}
\begin{proof}[Proof of Theorem \ref{Main Theorem 2}]  Let $m=n$ or $2n$ be as in
the Theorem \ref{Main Theorem 2} and $P^\chi(f)$ the Heegner point
construct in Theorem \ref{Main Theorem 2}.  By Theorem \ref{Main
Theorem 2},  we know that  $P^\chi(f)\in
E\left(\BQ(\sqrt{m^*})\right)^-\cong E^{(m)}(\BQ)$ is of infinite
order. Via the argument in the proof of Proposition \ref{Prop A},
the Euler system theory of Kolyvagin implies that the Mordell-Weil
group $E^{(m)}(\BQ)$ has rank one and the Shafarevich-Tate group
$\Sha(E^{(m)}/\BQ)$ is finite.  By the generalized Gross-Zagier
formula (\cite{YZZ} Theorem 1.2) and the non-vanishing of
$L(E^{(1)}, 1)$ and $L(E^{(2)}, 1)\neq 0$,  the vanishing order of
$L(s, E^{(m)})$ is exactly one.

It follows from (5.1) in Lemma 5.2 that
$$\rank_\BZ E^{(m)}(\BQ) +\dim_{\BF_2} \Sha(E^{(m)}/\BQ)[2]=1.$$
Since we have shown that $\rank_\BZ E^{(m)}(\BQ)=1$,
$\Sha(E^{(m)}/\BQ)[2]=0$, which implies that $\Sha(E^{(m)}/\BQ)$ has
odd cardinality.
\end{proof}

\end{document}